\documentclass[final,11pt]{article}
\usepackage{amsmath}
\usepackage{amsfonts}
\usepackage{makeidx}
\usepackage{amsthm}
\usepackage{mathtools}
\usepackage{amssymb}
\usepackage{bm}
\usepackage{cite}
\usepackage{showkeys}
\usepackage[english]{babel}
\usepackage{dblaccnt}
\usepackage{accents}
\usepackage{graphicx}
\usepackage{psfrag}
\usepackage{subfig}
\usepackage{color}
\usepackage{float}
\usepackage{amscd}
\usepackage[mathscr]{euscript}
\usepackage{lipsum}
\usepackage{epstopdf}
\usepackage{algorithm}
\usepackage{algpseudocode}

\newcommand{\R}{\ensuremath{\mathbb{R}}}

\newcommand{\x}{{\bf x}}
\newcommand{\y}{{\bf y}}
\renewcommand{\d}{\mathrm{d}}

\newtheorem{defi}{Definition}

\newtheorem{theorem}[defi]{Theorem}

\newtheorem{remark}[defi]{Remark}

\begin{document}

\title{Imaging of buried objects from multi-frequency experimental data
using a globally convergent inversion method}
\author{
Dinh-Liem Nguyen\thanks{Department of Mathematics and Statistics, University of North Carolina at
Charlotte, Charlotte, NC 28223, USA; (\texttt{dnguye70@uncc.edu}, \texttt{mklibanv@uncc.edu}, \texttt{lnguye50@uncc.edu})} 
\and Michael V. Klibanov\footnotemark[1] 
\and Loc H. Nguyen\footnotemark[1] 
\and Michael A. Fiddy\thanks{Optoelectronics Center, University of North Carolina at Charlotte,Charlotte, NC 28223, USA; 
(\texttt{mafiddy@uncc.edu})}  
}

\date{}
\maketitle

\begin{center}
\vspace{-0.6cm}
\textit{This paper is dedicated to Professor  Anatoly Borisovich
Bakushinsky, a distinguished expert in the theory of  ill-posed and inverse problems, on the occasion of his 80st birthday.}
\end{center}

\begin{abstract}
This paper is concerned with the numerical solution to a 3D coefficient 
inverse problem for buried objects with multi-frequency experimental data. 
The measured data, which are associated with a single direction of an incident 
plane wave, are backscatter data for targets buried in a sandbox.
These raw scattering data were collected using a microwave scattering facility at the University 
 of North Carolina at Charlotte.  
We develop a data preprocessing procedure and 
exploit a newly developed globally convergent inversion method for solving the 
inverse problem with these preprocessed data. 
 It is shown that dielectric constants of the buried targets as well as their
locations are reconstructed with a very good accuracy. We also prove a new analytical
result which rigorously justifies an important step of the so-called
``data propagation" procedure.
\end{abstract}

\sloppy


\bigskip


\section{Introduction}

\label{sec: intro}

Our target application is in imaging of explosive devices buried under
the ground, such as, e.g. antipersonnel land mines and improvised explosive
devices. Thus, we consider the model of a sandbox containing  objects which mimic
explosives. Our experimental data are the backscatter of the
electric field at different frequencies. Our goal
is to reconstruct locations and dielectric constants of scatterers from
these data using the solution of a Coefficient Inverse Problem (CIP). Since
the position of the source of the electric wave field is fixed, we have a
CIP with single measurement data. A new analytical result of this paper is
Theorem~\ref{thm propa} of Section~\ref{sect:SAR}, which rigorously justifies an
important step of the so-called  ``data propagation"
procedure used in our data preprocessing.

Of course, estimates of dielectric constants of targets do not
provide a sufficient information for separating explosives from clutters.
Still, these estimates might potentially be used as an additional piece of
information about explosive-like targets. Being combined with the currently
used energy information of radar images \cite{Soumekh:1999}, these estimates
might potentially lower the current false alarm rate. 

This work  is a continuation of the research of our group
about a globally convergent numerical method for solving a 
CIP~\cite{Kliba2016, Exp1, Exp2}. A CIP is the problem of reconstruction of at least one of the coefficients
of a partial differential equation (PDE) from the boundary measurements of its
solution. 
While the numerical method of~\cite{Kliba2016} was tested in~\cite{Exp1, Exp2}
on experimental data for targets placed in air, here we exploit it to solve 
a CIP for the case which is both more difficult and a more realistic one: 
this is the case of targets buried in a sandbox.
Compared with~\cite{Exp1, Exp2}, the obvious additional difficulty here is the
presence of the air/sand interface. Hence, the signal from the sand is mixed
with the signal from the target. The first thing, which comes in mind here,
is to consider a rather complicated mathematical model, which would take
this interface into account. However, our numerical method does not work for
this model, and it is yet unclear whether a corresponding modification can
ever be made. Hence, rather than considering this model, we preprocess the
experimental data in such a way that the preprocessed data are treated as
the ones for targets placed in air. We note that actually, there are various data
preprocessing procedures commonly used in engineering. 

As to the data preprocessing, we point out that verifying
performances of two versions of the globally convergent method for
experimental data~\cite{Klibanov:ip2010, Kuzhu2012, Kuzhuget:IEEE2013, Thanh2014, Exp1, Exp2}, 
it was established that it is
important to preprocess experimental data. The preprocessed data are then used
 as an input for the corresponding inversion algorithm. Although the data
preprocessing procedure is always a heuristic one, we believe that it is
justified by the accuracy of our reconstruction results, see Table 2. 

It is quite challenging to develop an efficient method for  solving CIPs numerically. 
This is because  CIPs are in general both highly nonlinear and ill-posed. Therefore, to solve
it numerically, it is important to address the following question: \emph{
How to obtain at least one point in a sufficiently small neighborhood of the
exact solution without any advanced knowledge of this neighborhood?} 
We call a numerical method for a CIP \emph{globally convergent} if
it addresses this question. More precisely, if: (1) a theorem is proved,
which positively addresses this question under a reasonable mathematical
assumption and (2) the error of the approximation of the true solution
depends only on the level of noise in the boundary data and on some
discretization errors, see~\cite{Beili2012, Beili2012a} for more details. 

Due to their real world applications, the CIPs have been studied
intensively. We refer the reader to several effective imaging methods that
can be used to detect the positions and the shapes of targets, such as, e.g.
level set methods, sampling methods, expansion methods, shape optimization
methods; see, e.g.~\cite{Ammar2004, Burge2005,
Cakon2006, Colto1996, Kirsc1998, LiLiuZou:smms2014,
Li2015} and references cited therein. The information about locations and shapes of targets
is not sufficient to identify  antipersonnel land mines and improvised
explosive devices. A more important piece of information is the value of the
dielectric constant of such a target, since this value characterizes a
material property of that target. Thus, we focus here on accurate
computations of values of dielectric constants of targets considered in our
experiments. As mentioned above, the knowledge of dielectric constants might 
serve as a piece of information for the radar community. To the best of our knowledge, 
this community relies only on the intensity of the radar images, see, 
e.g.~\cite{Kuzhu2012, Soumekh:1999},
and the differentiation between explosives and clutters is not yet
satisfactory.

A natural approach to solve a CIP is the optimization method, minimizing
some cost functionals, see, e.g. \cite{Bak,Chavent,Engl,Gonch,T}. This method is widely
used in the communities of inverse problems and engineering. However, those
cost functionals are in general not convex due to the nonlinearity of the CIP. 
This results in their multiple local
minima, see for example Figure 1 in \cite{Scales}. Consequently, a very good initial
guess for the true solution of the CIP is crucial to guarantee the
convergence of the minimizing sequence to that solution. A natural initial
guess for the desired coefficient is a constant function, assuming that the
true solution is close to the function  describing the background of the
medium. As a result, it is hard to reconstruct the targets with high
target/background contrasts. However, high contrasts are typical in
applications, see, e.g. the table of dielectric constants on the web site of
Honeywell with shortened link by Google https://goo.gl/kAxtzB, as well as
Table 2 of this paper. Motivated by this observation, we have developed in 
\cite{Kliba2016} a globally convergent algorithm which finds a good
approximation for the true solution of the CIP in the first stage and then
iteratively corrects such an \textquotedblleft initial guess" in the second
stage. Thus, the difficulty about the initial guess is resolved in the
algorithm of \cite{Kliba2016}. 


Our globally convergent numerical method works only for the case of
a CIP for the Helmholtz equation. It is yet unclear whether it can be
extended to a CIP for the Maxwell's system. Also, we have measured in our
experiments only a single component of the electric field. Indeed,
measurements of all three components would make our measurement system quite
more complicated. On the other hand, we want to have a simple measurement
system, which would be potentially used in a realistic environment of search
of explosives. Thus, we model the process of the propagation of the electric
field by the Helmholtz equation rather than by the full Maxwell's system. We
justify our model by the result of \cite{Beilina:cejm2013}, where it was
shown numerically that in the case of a rather simple medium (such as we
have) and the time dependent process the component of the electric field,
which was incident upon the medium, dominates two other components, and its
propagation is well governed by the time domain analog of the Helmholtz
equation. We also mention that in the above cited previous works on
experimental data~\cite{Klibanov:ip2010, Kuzhu2012, Kuzhuget:IEEE2013, 
Thanh2014, Exp1, Exp2}, either Helmholtz equation or its
analog for the case of the time domain was used to model the process.
Nevertheless, accurate results have been consistently obtained in all these
publications, so as in the current one. In our opinion, the latter provides
an additional justification of this modeling. Establishing a globally
convergent method for a CIP for the Maxwell's system is reserved for our
future research.

Recall that we work with the case of the backscattering data resulting from
a single measurement event generated by a single location of the source.
This means that we use only a minimal \textquotedblleft amount" of data. 
Except of the above cited works~\cite{Beili2012, Beili2012a, Klibanov:ip2010, 
Kuzhu2012, Kuzhuget:IEEE2013, Thanh2014, Kliba2016, Exp1, Exp2} and
corresponding references cited therein, the authors are unaware about other
globally convergent numerical methods for multidimensional CIPs, in which
the data are generated either by a single point source or by a single
direction of the incident plane wave. We refer the reader to some numerical
methods that solve the CIPs with the data generated by multiple directions
of the source, i.e. the Dirichlet to Neumann map, see \cite{Agalt2014,Kaban2004,Kab,Novikov}
and references cited therein. These methods work well without \textit{a
priori} knowledge of a small neighborhood of the exact coefficient.

The paper is organized as follows. In the next two sections, we formulate both our 
CIP and the globally convergent method developed in \cite{Kliba2016}. Section 4 is devoted 
to our justification for the data propagation process. We next describe the
experimental set up and the data preprocessing procedure in Section 
\ref{sect:imag} and present the numerical results in Section \ref{sect:result}.

\section{The problem formulation}

\label{sec: 2}

Let $\Omega $ be an open and bounded domain in $\mathbb{R}^{3}$ with its smooth
boundary $\partial \Omega$. For $\mathbf{x}=(x,y,z)^\top\in \mathbb{R}^{3}$, 
let $c(\mathbf{x})$ be the spatially distributed dielectric constant of the medium. Assume that 
$c(\mathbf{x})$ is a smooth function and that
\begin{equation}
c(\mathbf{x}) = 1 \text{ for }\mathbf{x}\in \mathbb{R}^{3}\setminus \overline{\Omega } \quad  \text{and } 1 \leq c(\mathbf{x}) \leq C, 
\end{equation}
where  $C$ is a positive constant. Denoting by  $k>0$ the wavenumber, we consider the incident plane wave 
\begin{equation}
u_{\mathrm{inc}}(\mathbf{x},k)=\exp (ikz).  \label{eqn incident}
\end{equation}%
The total wave $u(\mathbf{x},k)$ is governed by the Helmholtz
equation 
\begin{equation}
\Delta u(\mathbf{x},k)+k^{2}c(\mathbf{x})u(\mathbf{x},k)=0,\quad \mathbf{x}
\in \mathbb{R}^{3}  \label{eqn Hel}
\end{equation}%
with the Sommerfeld radiation condition 
\begin{equation}
\partial _{r}u_{\mathrm{sc}}(\mathbf{x},k)-iku_{\mathrm{sc}}(\mathbf{x},k)=o(r^{-1}),\quad r=|\mathbf{x}|\rightarrow \infty ,
\label{eqn outgoing}
\end{equation}%
where $u_{\mathrm{sc}}(\mathbf{x},k)$ denotes the scattered wave
\begin{equation}
u_{\mathrm{sc}}(\mathbf{x},k)=u(\mathbf{x},k)-u_{\mathrm{inc}}(\mathbf{x},k).
\label{eqn usc}
\end{equation} 
Let $\Gamma \subset \partial \Omega $ be the  ``backscatter" part of $\Omega $ (see Figure \ref{fi:1}(b) for an illustration). 
We are interested in the following coefficient inverse problem:

\noindent \textbf{Coefficient Inverse Problem (CIP).} Let $
\underline{k}$ and $\overline{k}$ be two numbers such that $0<
\underline{k}<\overline{k}.$ Determine $c(\mathbf{x})$ for $\mathbf{x}\in
\Omega $ from the boundary measurement 
\begin{equation}
g(\mathbf{x},k) :=u(\mathbf{x},k), \quad \mathbf{x}\in \Gamma,
k\in [ \underline{k},\overline{k}],
\end{equation}
where $u(\mathbf{x},k)$ is the total wave associated to the 
incident plane wave $u_{\mathrm{inc}}(\mathbf{x},k)=\exp (ikz)$.

We note that the multi-frequency data in the CIP correspond to
a single direction of the incident plane wave. Therefore, we call this problem
a CIP with single measurement. 


\begin{remark}
The function $g(\mathbf{x},k)$ is given on the backscatter side 
$\Gamma $ of $\partial \Omega $. However, the theory of \cite{Kliba2016} works only for the case when
measurements are conducted on the entire boundary $\partial \Omega$. 
Hence, in the numerical implementation the
data are heuristically complemented on the entire boundary $\partial \Omega$
(see more details in section~\ref{sect:completion}.)
%
%
%
\end{remark}

It is well-known that the forward problem (\ref{eqn incident})--(\ref{eqn usc}) can be
reformulated as the Lippmann-Schwinger integral equation. Since this integral equation 
plays an important role in both the analysis and the numerical implementation 
of our inversion method developed in~\cite{Kliba2016}, it is briefly presented
in the following. Assume that the spatially distributed dielectric
constant $c(\mathbf{x}) $ is known. It follows from \eqref{eqn
Hel} and \eqref{eqn
usc} that 
\begin{equation}
\Delta u_{\mathrm{sc}}(\mathbf{x},k)+k^{2}u_{\mathrm{sc}}(\mathbf{x},k)=-k^{2}\beta (\mathbf{x})u(\mathbf{x},k),  \label{7}
\end{equation}%
where the function 
\begin{equation}
\beta (\mathbf{x})=c(\mathbf{x})-1,\quad \mathbf{x}\in \mathbb{R}^{3}
\end{equation}%
is supported in $\Omega$. Equation \eqref{7} and the  radiation
condition \eqref{eqn outgoing} deduce that the function $u_{\mathrm{sc}}$
satisfies 
\begin{equation}
u_{\mathrm{sc}}(\mathbf{x},k)=k^{2}\int_{\mathbb{R}^{3}}\Phi_k(\mathbf{x},
\mathbf{y})\beta (\mathbf{y})u(\mathbf{y},k)\d\mathbf{y,}\quad \mathbf{x}%
\in \mathbb{R}^{3}.  \label{10}
\end{equation}%
where 
\begin{equation}
\Phi_k (\mathbf{x},\y)=\frac{\exp (ik|\mathbf{x}-\y|)}{4\pi |\mathbf{x}-\y|},
\quad \mathbf{x} \neq \y,
\end{equation}%
is the fundamental solution of the Helmholtz equation in $3$D. From 
\eqref{eqn usc} and \eqref{10}, we have the Lippmann-Schwinger equation 
\begin{equation}
u(\mathbf{x},k)=u_{\mathrm{inc}}(\mathbf{x},k)+k^{2}\int_{\Omega}
\Phi_k(\mathbf{x},\mathbf{y})\beta (\mathbf{y})u(\mathbf{y},k)\d\mathbf{y,}\quad 
\mathbf{x}\in \mathbb{R}^{3}.  \label{eqn LS}
\end{equation}%
The integral in \eqref{eqn LS} is written on the domain $\Omega $
because supp$(\beta )\subset \Omega .$ The above derivation of the
Lippmann-Schwinger equation \eqref{eqn LS} is  from \cite[Chapter 8]{Colto1998}. 
We also refer the reader to \cite{Colto1998} for the well-posedness  of \eqref{eqn LS}. 
The forward problem corresponding to our CIP is equivalent to equation \eqref{eqn LS} for each value of
the wavenumber $k\in \lbrack \underline{k},\overline{k}]$.

\section{The globally convergent numerical method}

\label{sect:FP4} 

For the convenience of the reader, we briefly summarize in this section the
method of \cite{Kliba2016}, which we will use in this paper with
appropriate modifications. In this section, we assume that the total wave
field 
\begin{equation}
u(\mathbf{x},k)\neq 0, \quad \forall \mathbf{x}\in \overline{\Omega },\forall k\in
\lbrack \underline{k},\overline{k}].  \label{11}
\end{equation}
This assumption is justified in \cite{Kliba2016} for a sufficiently large
number $\underline{k}$. We also observed from our numerical study 
for both simulated and experimental data that (\ref{11}) is also true
for all wavenumbers. 

\subsection{An integro-differential equation}

\label{sect:integral}

In this subsection, we establish a boundary value problem for an
integro-differential equation. The solution of this problem directly
provides the desired solution of the above CIP. Since the vector function $%
\nabla u(\mathbf{x},k)/u(\mathbf{x},k)$, $\mathbf{x}\in \Omega $, is curl
free for each $k\in \lbrack \underline{k},\overline{k}]$, one can explicitly
find a smooth function $v(\mathbf{x},k)$ such that (see \cite{Kliba2016}) 
\begin{equation}
u(\mathbf{x},k)=\exp (v(\mathbf{x},k)),\quad \mathbf{x}\in \Omega ,k\in
\lbrack \underline{k},\overline{k}].  \label{12}
\end{equation}%
Plugging (\ref{12}) into \eqref{eqn Hel} gives 
\begin{equation}
\Delta v(\mathbf{x},k)+(\nabla v(\mathbf{x},k))^{2}=-k^{2}c(\mathbf{x}).
\label{eqn v}
\end{equation}%
In order to eliminate $c(\mathbf{x})$ from (\ref{eqn v}), we differentiate 
\eqref{eqn v} with respect to $k$ and get 
\begin{equation}
\Delta q(\mathbf{x},k)+2\nabla q(\mathbf{x},k)\cdot v(\mathbf{x},k)=\frac{%
2\Delta v(\mathbf{x},k)+2(\nabla v(\mathbf{x},k))^{2}}{k},\quad \mathbf{x}%
\in \Omega ,k\in \lbrack \underline{k},\overline{k}],  \label{14}
\end{equation}%
where the function $q(\mathbf{x},k)$ is defined as 
\begin{equation}
q(\mathbf{x},k)=\partial _{k}v(\mathbf{x},k),\quad \mathbf{x}\in \Omega
,k\in \lbrack \underline{k},\overline{k}].  \label{15}
\end{equation}%
Denoting by $V(\mathbf{x})$ the ``tail function" 
\begin{equation}
V(\mathbf{x})=v(\mathbf{x},\overline{k}),
\end{equation}%
we have
\begin{equation}
v(\mathbf{x},k)=-\int_{k}^{\overline{k}}q(\mathbf{x},s)\d s+V(\mathbf{x}%
),\quad \mathbf{x}\in \Omega ,k\in \lbrack \underline{k},\overline{k}].
\label{17}
\end{equation}%
Plugging \eqref{17} into \eqref{14}, we obtain 
\begin{equation*}
\frac{k}{2}\Delta q(\mathbf{x},k)+k\nabla q(\mathbf{x},k)\cdot \left(
-\int_{k}^{\overline{k}}\nabla q(\mathbf{x},s)\d s+\nabla V(\mathbf{x})\right) 
\end{equation*}%
\begin{equation}
=-\int_{k}^{\overline{k}}\Delta q(\mathbf{x},s)\d s+\Delta V(\mathbf{x}%
)+\left( -\int_{k}^{\overline{k}}\nabla q(\mathbf{x},s)\d s+\nabla V(\mathbf{x}%
)\right) ^{2},  \label{eqn ide}
\end{equation}%
which is the integro-differential equation we want to derive in this
subsection. The function $q$ also satisfies the boundary condition 
\begin{equation}
q(\mathbf{x},k)=\psi (\mathbf{x},k)=\frac{\partial _{k}g(\mathbf{x},k)}{g(%
\mathbf{x},k)},\quad \mathbf{x}\in \partial \Omega .  \label{20}
\end{equation}

In the next subsection, we discuss an approximation for the vector function $%
\nabla V(\mathbf{x})$.

\subsection{The initial approximation $V_{0}( \mathbf{x}) $ for
the tail function}

\label{sec: tail}

The vector function $\nabla V(\mathbf{x})$ is an important
component of equation \eqref{eqn ide}. Therefore, approximating it is
crucial in our inversion algorithm. In this subsection, following \cite{Kliba2016}, we present the
initial approximation $\nabla V_{0}( \mathbf{x}) $ for 
$\nabla V(\mathbf{x}).$ 
The follow up approximations $\nabla
V_{n}( \mathbf{x}) $ are computed in our iterative process.
For the concept of the global convergence mentioned in the 
Introduction, it is
crucial to find $\nabla V_{0}(\mathbf{x})$, without using any 
advanced knowledge of a small neighborhood of the exact
solution of our CIP.

We approximate the function $\nabla V( \mathbf{x})$ using the
following asymptotic expansion of the function $u(\mathbf{x},k)$
with respect to $k$  
\begin{equation}
u(\mathbf{x},k)=A(\mathbf{x})\exp (ik\tau (\mathbf{x}))(1+O(1/k)),\quad \mbox{as }k\rightarrow \infty, \mathbf{x}\in \Omega 
\end{equation}%
for some positive functions $A(\mathbf{x})$ and $\tau(\mathbf{x})$, see
details in \cite{Kliba2016}. This expansion is uniform for all $\mathbf{x}%
\in \Omega .$ 

Recall that the tail function is defined as $V(\mathbf{x})=v(\mathbf{x},%
\overline{k})$. Assuming that $\underline{k}$ and $\overline{k}$ are sufficiently large, 
using \eqref{12} and dropping the term $O(1/k)$, we obtain 
\begin{equation}
V(\mathbf{x})\approx \ln A(\mathbf{x})+i\overline{k}\tau (\mathbf{x}%
)\approx i\overline{k}\tau (\mathbf{x}), \quad \mathbf{x}\in \Omega.
\end{equation}%
Therefore, the function $q(\mathbf{x},\overline{k})$ can be approximated as 
\begin{equation}
q(\mathbf{x},\overline{k})=\partial _{k}v(\mathbf{x},k)\mid_{k=\overline{k}}
\approx i\tau (\mathbf{x})\approx V(\mathbf{x})/\overline{k}, \quad \mathbf{x}\in \Omega.
\label{21}
\end{equation}%
Plugging \eqref{21} into \eqref{eqn ide} and setting in \eqref{eqn ide} $k=%
\overline{k}$, we obtain 
\begin{equation}
\frac{1}{2}\Delta V(\mathbf{x})+(\nabla V(\mathbf{x}))^{2}=\Delta V(\mathbf{x})+\left( \nabla V(\mathbf{x})\right) ^{2},  \quad \mathbf{x}\in \Omega 
\end{equation}%
which yields 
\begin{equation}
\Delta V(\mathbf{x})=0, \quad \mathbf{x}\in \Omega .
\end{equation}
Since only $\nabla V(\mathbf{x})$ and $\Delta V(\mathbf{x})=\mathrm{div}(\nabla V(\mathbf{x}))$ are involved in equation 
(\ref{eqn ide}), instead of finding $V(\mathbf{x})$ we compute directly the vector function $\nabla V(\mathbf{x})$ in our numerical implementation 
 by solving the following problem
\begin{equation}
\label{tail}
\left\{ 
\begin{array}{rcll}
\Delta (\nabla V(\mathbf{x}))& = & 0 & \mbox{in }\Omega,  \\ 
\nabla V(\mathbf{x}) & = & R(\x,\overline{k}) & \mbox{on }\partial \Omega,
\end{array}
\right. 
\end{equation}
where  $R(\x,\overline{k})$ is a certain vector function, which is known approximately. 
We refer to~\cite[Section 7.4]{Kliba2016}
for all the details about approximation of $R(\x,\overline{k})$ on the entire boundary $\partial \Omega$.
We consider the solution of  problem~\eqref{tail} as the first
approximation $\nabla V_{0}$ of the
vector function $\nabla V$.

\subsection{The globally convergent algorithm}

\label{sect:Algorithm1}

In this subsection, we describe our inversion method for the
CIP. For $N \in  \mathbb{N}$, consider the uniform partition 
\begin{equation}
k_{N}=\underline{k}<k_{N-1}<\dots <k_{1}<k_{0}=\overline{k}  \label{100}
\end{equation}
of the interval $[\underline{k},\overline{k}]$ with the step size $%
h=k_{i-1}-k_{i}$, $1\leq i\leq N.$ For each $n\in \{1,\dots ,N\},$ denote 
\begin{equation}
q_{n}(\mathbf{x})=q(\mathbf{x},k_{n}),\quad u_{n}(\mathbf{x})=u(\mathbf{x},k_{n}), \quad \mathbf{x}\in \Omega.  \label{101}
\end{equation}
Recall that the first approximation $\nabla V_{0}$ for the gradient $\nabla V
$ of the tail function is constructed in Section \ref{sec: tail}. We
assume that $\nabla V_{n-1}$ is known, which implies that $\Delta V_{n-1}=%
\mathrm{div}\left( \nabla V_{n-1}\right) $ is known as well, where $n\in
\{1,\dots ,N\}$. By (\ref{100}) and (\ref{101}) the discrete form of
equation \eqref{eqn ide} is 
\begin{equation*}
\frac{k_{n}}{2}\Delta q_{n}(\mathbf{x})+k_{n}\nabla q_{n}(\mathbf{x})\cdot
\left( -\nabla Q_{n-1}(\mathbf{x})+\nabla V_{n-1}(\mathbf{x})\right) 
\end{equation*}%
\begin{equation}
=-\Delta Q_{n-1}(\mathbf{x})+\Delta V_{n-1}(\mathbf{x})+\left( -\nabla
Q_{n-1}(\mathbf{x})+\nabla V_{n-1}(\mathbf{x})\right) ^{2},  \quad \mathbf{x}\in \Omega, \label{30}
\end{equation}%
where 
\begin{equation}
Q_{n-1}(\mathbf{x})=h\sum_{i=0}^{n-1}q_{n}(\mathbf{x}),\quad \mathbf{x}\in
\Omega .  \label{31}
\end{equation}%
Here, we approximate the integral 
\begin{equation*}
\int_{k}^{\overline{k}}q(\mathbf{x},s)\d s
\end{equation*}
by $Q_{n-1}(\mathbf{x})$ instead of $Q_{n}(\mathbf{x})$ to remove the
nonlinearity of \eqref{eqn ide}. The resulting error is $O(h),$ as $%
h\rightarrow 0.$ The boundary condition for the function $q_{n}(\mathbf{x})$
is 
\begin{equation}
q_{n}(\mathbf{x})=\frac{g(\mathbf{x},k_{n})-g(\mathbf{x},k_{n+1})}{hg(%
\mathbf{x},k_{n})},\quad \mathbf{x}\in \partial \Omega .  \label{32}
\end{equation}

\begin{remark}
 Even though the derivative of the data with respect
to $k$ is calculated in (\ref{32}) via the finite difference, we
have not observed any instabilities in our computations, probably because
the step size $h$ was not exceedingly small. The same can be stated
about all above cited publications about the globally convergent methods of
this group.
\end{remark}

In short, the algorithm is based on the following iterative process: (1)
given $\nabla V_{n-1},\Delta V_{n-1},$ solve the  Dirichlet boundary
value problem (\ref{30})--(\ref{32}); (2) use (\ref{eqn v}) to calculate the
function $c_{n}(\mathbf{x}) $ via $q_{n}(\mathbf{x}), \nabla V_{n-1}(\mathbf{x})$ 
and $\Delta V_{n-1}(\mathbf{x}) $; (3) update $\nabla V_{n}(\mathbf{x}) $ and $%
\Delta V_{n}(\mathbf{x}) $ by solving the forward problem 
(\ref{eqn LS}) with $\beta(\mathbf{x}) :=c_{n}(\mathbf{x}) -1$. To increase the stability of this iterative process, we
arrange  internal iterations inside steps (2) and (3). The whole algorithm is
summarized as follows, see \cite{Kliba2016, Exp1, Exp2} for more details. 
\begin{algorithm}[H]
\caption{Globally convergent algorithm}
\label{alg:globalconv}
\begin{algorithmic}[1]
\State{Given $\nabla V_{0}$, set $q_{0}:=0$}
\For{$n = 1, 2, \dots, N$}
\State{Set $q_{n, 0} := q_{n - 1}$ and $\nabla V_{n, 0} := \nabla V_{n - 1}$}
\For{$i=1, 2, \dots, I_N$}
\State{Find $q_{n,i}$ by solving the boundary value problem~\eqref{30}--\eqref{32}}.
\State{Update 
$\nabla v_{n,i} :=- ( h\nabla q_{n,i}
+ h\sum_{j=0}^{n-1}\nabla q_{j}) +\nabla V_{n,i-1}$ in $\Omega$}.
\State{Update $c_{n,i}$ via~\eqref{eqn v}}.
\State{Find $u_{n,i}(\x,\overline{k})$ by solving integral equation~\eqref{eqn LS} in $\Omega$ with $\beta(\mathbf{x}) :=c_{n,i}( \mathbf{x}) -1.$}
\State{Update  $\nabla V_{n,i}(\x) := \nabla u_{n,i}(\x,\overline{k})/ u_{n,i}(\x,\overline{k})$}.
\EndFor
\State{Update $q_{n} :=q_{n,I_N}$, $c_n :=c_{n,I_N}$ and $\nabla V_n := \nabla V_{n,I_N}$}.
\EndFor
\State{Choose $c$ by \emph{the-criterion-of-choice}}
\end{algorithmic}
\end{algorithm}

\begin{remark}
In Algorithm~\ref{alg:globalconv}
``the-criterion-of-choice" means the stopping rule and 
choosing the final result presented in Section~\ref{sect:stopping}.
This rule was computationally developed in~\cite{Kliba2016, Exp2}. 
 Theorem 6.1 of \cite{Kliba2016} ensures the
global convergence of this algorithm. 
\end{remark}

\section{A justification of the data propagation formula}

\label{sect:SAR} 

This section is devoted to a rigorous mathematical justification of
one step of the data propagation procedure, which is used both in the
current paper and in~\cite{Kliba2016, Exp1, Exp2}. This procedure is also known in the optics
community under the name  ``the angular spectrum
representation", see, e.g.~\cite[Chapter 2]{Novot2012}. However, this
important step was used only heuristically so far and was never rigorously
justified. For the convenience of the presentation 
in this section, we use the notations:
$$
\mathbf{x} = (x_1,x_2,x_3)^\top, \quad \mathbf{y} = (y_1,y_2,y_3)^\top.
$$
In addition, for the analysis in this section, it is important to note that the validity of all the integral transforms,
the convolution and their related properties can be understood in the sense of 
tempered distributions, see, e.g.~\cite[Chapter 2]{Vladi1979}. The integrands in this section
can be verified as tempered distributions. Indeed, for example, if we consider function 
$\phi(\x) =\exp(ik|\x|)/|\x|$, when $x_3$ is away from 0 which is 
the case of our analysis below, then 
we have
$$
\int_{\R^2} \frac{|\phi(x_1,x_2,x_3)|}{(1+ \sqrt{x_1^2+x_2^2})^3} \d x_1 \d x_2 < \infty.
$$
The latter fact defines $\phi$ as a bounded linear functional $\phi(\psi) = \int_{\R^2} \phi(\x) \psi(x_1,x_2) \d x_1\d x_2$
on the space of rapidly decreasing functions on $\R^2$. We recall that the scattered field in our problem behaves
like function $\phi(\x)$ on some (measurement) plane away from the scatterer. Other integrands can be verified 
as tempered distributions similarly.

\subsection{Data propagation formula}

\label{sect:4.1}

We now describe our the data propagation process.  Suppose that the domain $\Omega
\subset \left\{ x_3>0\right\}$. Hence, the half space $\left\{ x_3<0\right\}\subset \mathbb{R}^{3}$ 
is homogeneous, i.e. $c(\mathbf{x}) =1$
for $\mathbf{x}\in \left\{ x_3<0\right\}$. Then it follows from  (\ref{eqn
incident})--(\ref{eqn usc}) that the function $u_{sc}$ is actually a 
backscattered wave in $\left\{ x_3<0\right\} $ and it satisfies the following
conditions
\begin{equation}
\Delta u_{\mathrm{sc}}(\mathbf{x},k)+k^{2}u_{\mathrm{sc}}(\mathbf{x},k) = 0,\quad \mathbf{x}\in \{ x_3<0\},  \label{4.1}
\end{equation}
\begin{equation}
\partial _{r}u_{\mathrm{sc}}(\mathbf{x},k)-iku_{\mathrm{sc}}(%
\mathbf{x},k)=o(r^{-1}),\quad r=|\mathbf{x}|\rightarrow \infty ,  \label{4.2}
\end{equation}%
Let%
\begin{equation}
u_{\mathrm{sc}}(x_1,x_2,0,k)=U( x_1,x_2,k).  \label{4.3}
\end{equation}%
Let $d<0$ be a number. Suppose that we measure the function 
$u_{\mathrm{sc}}(\mathbf{x},k)$ at the plane $\left\{ x_3 = d\right\}$,
\begin{equation}
u_{\mathrm{sc}}(x_1,x_2,d,k)=V( x_1,x_2,k).  \label{4.4}
\end{equation}
The goal of the data propagation procedure is to ``move" the
data (\ref{4.4}) closer to the domain of interest $\Omega .$ In other words,
given the function $V$, we want to approximate the function $U$ in (\ref{4.3}).

Consider the Fourier transform of the scattered field with respect to $x_1, x_2$ 
\begin{equation}
\widehat{u}_{\mathrm{sc}}(v_{1},v_{2},x_3,k)= \int_{\mathbb{R}^{2}}u_{\mathrm{sc%
}}(\mathbf{x},k)e^{-i(x_1v_{1}+x_2v_{2})} \d x_1 \d x_2.  \label{4.5}
\end{equation}%
Then (\ref{4.1}) and (\ref{4.3}) imply that
\begin{equation}
\partial_{x_3}^{2}\widehat{u}_{\mathrm{sc}}+(
k^{2}-v_{1}^{2}-v_{2}^{2}) \widehat{u}_{\mathrm{sc}}=0, \quad x_3<0,
\label{4.6}
\end{equation}%
\begin{equation}
\widehat{u}_{\mathrm{sc}}(v_{1},v_{2},0,k)=\widehat{U }(v_{1},v_{2},k).
\label{4.7}
\end{equation}%
Also, by (\ref{4.4})%
\begin{equation}
\widehat{u}_{\mathrm{sc}}(v_{1},v_{2},d,k)=\widehat{V }(v_{1},v_{2},k).
\label{4.8}
\end{equation}%
Here $\widehat{U}$ and $\widehat{V}$ are Fourier transforms (\ref{4.5}) of
functions $U $ and $V $ respectively.
We now solve the problem (\ref{4.6})--(\ref{4.7}). Consider two cases:

\textbf{Case 1}. $k^{2}-v_{1}^{2}-v_{2}^{2}<0.$ Then (\ref{4.6}) implies
that 
\begin{equation}
\widehat{u}_{\mathrm{sc}}(v_{1},v_{2},x_3,k)=C_{1}e^{-\sqrt{%
v_{1}^{2}+v_{2}^{2}-k^{2}}x_3}+C_{2}e^{\sqrt{v_{1}^{2}+v_{2}^{2}-k^{2}}x_3}, \quad x_3<0,
\label{4.9}
\end{equation}%
where $C_{1}$ and $C_{2}$ are independent on $z$.

\textbf{Case 2}. $k^{2}-v_{1}^{2}-v_{2}^{2}>0.$ Then (\ref{4.6}) implies that%
\begin{equation}
\widehat{u}_{\mathrm{sc}}(v_{1},v_{2},x_3,k) = 
C_{1}e^{-i\sqrt{k^{2}-v_{1}^{2}-v_{2}^{2}}x_3}+C_{2}e^{i\sqrt{k^{2}-v_{1}^{2}-v_{2}^{2}}x_3}, \quad x_3<0.
\label{4.10}
\end{equation}

The problem with formulae (\ref{4.9}) and (\ref{4.10}) is that it is unclear
how to connect $C_{1}$ and $C_{2}$ with the function $\widehat{U}(v_{1},v_{2},k)$ in (\ref{4.7}). 
Indeed, it is unclear how to use the
radiation condition (\ref{4.2}) after the Fourier transform (\ref{4.5}). In
the past works~\cite{Exp1,Exp2} it
was set in (\ref{4.9}) and (\ref{4.10}) that
\begin{equation}
\widehat{u}_{\mathrm{sc}}(v_{1},v_{2},x_3,k)=\left\{ 
\begin{array}{c}
\widehat{U }(v_{1},v_{2},k)e^{\sqrt{v_{1}^{2}+v_{2}^{2}-k^{2}}x_3}, \quad x_3<0
\text{ in Case 1,} \\ 
\widehat{U}(v_{1},v_{2},k)e^{-i\sqrt{k^{2}-v_{1}^{2}-v_{2}^{2}}x_3}, \quad x_3<0
\text{ in Case 2.}%
\end{array}%
\right.   \label{4.11}
\end{equation}%
However, it was unclear until now how to justify (\ref{4.11}) rigorously.
Theorem~\ref{thm propa} in Section~\ref{sect:4.2} provides the desired rigorous justification of
(\ref{4.11}).

In any case, assume that (\ref{4.11}) is true. Then in Case 1 the function $
\widehat{u}_{\mathrm{sc}}(v_{1},v_{2},x_3,k)$ decays exponentially with
respect to $x_3$\ thus becoming sufficiently small at $x_3=d$, if $|d|$ is
sufficiently large. Hence, just as in \cite[Chapter 2]{Novot2012}, we ignore
Case 1. Next, considering only Case 2 and using (\ref{4.8}), we obtain
\begin{equation*}
\widehat{U}(v_{1},v_{2},k)=\widehat{V}(v_{1},v_{2},k)e^{i\sqrt{k^{2}-v_{1}^{2}-v_{2}^{2}}d}.
\end{equation*}%
Hence, 
\begin{equation}
U (x_1,x_2,k)=\frac{1}{2\pi }\int\limits_{v_{1}^{2}+v_{2}^{2}<k^{2}}%
\widehat{V }(v_{1},v_{2},k)e^{-i(x_1v_{1}+x_2v_{2})}e^{i\sqrt{k^{2}-v_{1}^{2}-v_{2}^{2}}
d}\d v_{1}\d v_{2}.  \label{4.12}
\end{equation}%
Formula (\ref{4.12}) is the desired formula of the data propagation
procedure.

\subsection{Rigorous justification of (\protect\ref{4.11})}

\label{sect:4.2}


The following theorem rigorously justifies formula (\ref{4.11}). 

\begin{theorem}
\label{thm propa} Suppose that  $\varphi (x_1,x_2,k)$ is a bounded measurable function 
with respect to $x_1$ and $x_2$, which has compact support. Let the
function $w( \mathbf{x},k) $ satisfy the following conditions
\begin{align}
 \label{eq:Helmeq}
\Delta w(\mathbf{x},k)+k^{2}w(\mathbf{x},k) & =0\quad \text{ in }\{x_3<0\},  \\
w(\mathbf{x},k)& =\varphi \quad \text{ on }\{x_3=0\}, \\
 \label{eq:Helmeq3}
 \partial_r w(\mathbf{x},k) -ikw(\mathbf{x},k) & =o(r^{-1}),\quad r=|%
\mathbf{x}|.
\end{align}%
Then the Fourier transform of $w$
with respect to $x_1$ and $x_2$, that is, 
\begin{equation*}
\widehat{w}(v_{1},v_{2},x_3,k) =  \int_{\mathbb{R}^{2}}w(\mathbf{x},k)e^{-i(x_1v_{1}+x_2v_{2})}\d x_1\d x_2
\end{equation*}%
is given by 
\begin{equation}
\widehat{w}(v_{1},v_{2},x_3,k)=\widehat{\varphi }(v_{1},v_{2},k)e^{-i\kappa x_3}\quad \text{for all }x_3 < 0,  \label{eq:uhat}
\end{equation}%
where 
\begin{equation*}
\kappa =%
\begin{cases}
\sqrt{k^{2}-v_{1}^{2}-v_{2}^{2}}, & v_{1}^{2}+v_{2}^{2}\leq k^{2}, \\ 
i\sqrt{v_{1}^{2}+v_{2}^{2}-k^{2}}, & v_{1}^{2}+v_{2}^{2}>k^{2}.
\end{cases}%
\end{equation*}
\end{theorem}

\begin{remark}
The assumption about the compact support of $\varphi$ is simply a technical assumption
in our proof.  However, it can be also justified by the fact that our measured data are only available 
on a rectangle in the $x_1x_2$-plane. 
We note that $\widehat{w}(v_{1},v_{2},x_3,k)$ in~\eqref{eq:uhat} satisfies the
radiation condition 
\begin{equation*}
\partial_{x_3}\widehat{w}+i\kappa \widehat{w}=0,
\end{equation*}%
which also means that it is a plane wave propagating to the negative $x_3-$direction.
\end{remark}

\begin{proof}
For $ \x,\y \in \R^3$, recall the fundamental solution of problem~\eqref{eq:Helmeq}--\eqref{eq:Helmeq3}
$$
\Phi_k(\x,\y) = \frac{e^{ik|\x-\y|}}{4\pi|\x-\y|}, \quad \x \neq \y.
$$
Then the fundamental solution $G$ of the Helmholtz equation  in  half-space $\{y_3<0\}$ is given by
\begin{align*}
G(\x,\y) = \Phi_k(\x,\y) - \Phi_k(\x,\y'),
\end{align*}
where $\y = (y_1,y_2,y_3)$ and $\y' = (y_1,y_2,-y_3)$. 
The solution $w$ in~\eqref{eq:Helmeq}--\eqref{eq:Helmeq3} can be written as the double layer potential 
$$
w(\x) = -\int_{\R^2} \frac{\partial G(\x,\y)}{\partial y_3}\left|_{\{y_3=0\}}\right. \varphi(y_1,y_2) \d y_1 \d y_2.
$$
The integral above is a convolution of a tempered distribution and a distribution with compact support.
This convolution is well defined and we have the following (see~\cite[Chapter 2]{Vladi1979})
\begin{align*}
&\widehat{w}(v_1,v_2,x_3) = -\frac{\widehat{\partial G}}{\partial y_3}(v_1,v_2,x_3) \,\widehat{\varphi}(v_1,v_2) \\
&=-\widehat{\varphi}(v_1,v_2) \int_{\R^2} \frac{\partial}{\partial y_3} \big[ \Phi(z_1,z_2,x_3-y_3) - \Phi(z_1,z_2,x_3+y_3) \big]\left|_{\{y_3=0\}}\right. e^{-i(z_1v_1+z_2v_2)} \d z_1 \d z_2 \\
&=-\widehat{\varphi}(v_1,v_2)\frac{\partial}{\partial y_3} \int_{\R^2}  \big[ \Phi(z_1,z_2,x_3-y_3) - \Phi(z_1,z_2,x_3+y_3) \big] e^{-i(z_1v_1+z_2v_2)} \d z_1 \d z_2
\left|_{\{y_3=0\}}\right..
\end{align*}
Now we compute 
\begin{align*}
I_1 &= \int_{\R^2} \Phi(z_1,z_2,x_3-y_3)e^{-i(z_1v_1+z_2v_2)} \d z_1 \d z_2 \\ 
& =\int_\R\int_\R \frac{e^{ik\sqrt{z_1^2+z_2^2+(x_3-y_3)^2}}}{4\pi\sqrt{z_1^2+z_2^2+(x_3-y_3)^2}} e^{-iz_1v_1} \d z_1 e^{-iz_2v_2} \d z_2 \\
&=\int_\R \int_\R \frac{e^{ik\sqrt{z_1^2+z_2^2+(x_3-y_3)^2}}}{4\pi\sqrt{z_1^2+z_2^2+(x_3-y_3)^2}} \cos(z_1v_1) \d z_1 e^{-iz_2v_2} \d z_2.
\end{align*}
From the table of integral transforms~\cite{Batem1954}  we have 
\begin{align*}
 &\int_0^\infty \frac{e^{ik\sqrt{z_1^2+z_2^2+(x_3-y_3)^2}}}{\sqrt{z_1^2+z_2^2+(x_3-y_3)^2}} \cos(z_1v_1) \d z_1  \\
&=  
\begin{dcases}
-\frac{\pi}{2} Y_0\big[(z_2^2+(x_3-y_3)^2)^{1/2}(k^2-v_1^2)^{1/2} \big]+ i \frac{\pi}{2}J_0\big[(x_2^2+(x_3-y_3)^2)^{1/2}(k^2-v_1^2)^{1/2} \big] \\
K_0\big[(z_2^2+(x_3-y_3)^2)^{1/2}(v_1^2- k^2)^{1/2} \big]
\end{dcases} \\
 &= 
 \begin{dcases}
 \frac{i\pi}{2} H^{(1)}_0\big[(z_2^2+(x_3-y_3)^2)^{1/2} (k^2-v_1^2)^{1/2} \big], \quad  v_1^2<k^2 \\
K_0\big[(z_2^2+(x_3-y_3)^2)^{1/2}(v_1^2- k^2)^{1/2} \big], \quad  v_1^2>k^2.
\end{dcases} 
\end{align*}
Again from~\cite{Batem1954} and the following formulae from~\cite{Grads2007} 
$$ H^{(1)}_{-1/2} (s) =  \sqrt{\frac{2}{\pi s}} e^{is} , \quad K_{-1/2}(s) = \sqrt{\frac{\pi}{2s}} e^{-s} $$
we have for $y_3 >x_3$ that
\begin{align*}
 &\int_0^\infty  \frac{i\pi}{2}H^{(1)}_0\big[(z_2^2+(x_3-y_3)^2)^{1/2}(k^2-v_1^2)^{1/2} \big] \cos(z_2v_2) \d z_2  \\
& = \begin{dcases}
             \frac{i\pi e^{-i(x_3-y_3)\sqrt{k^2-v_1^2-v_2^2}}}{2\sqrt{k^2-v_1^2-v_2^2}}, & v_1^2+v_2^2 <k^2, \\
            \frac{\pi e^{(x_3-y_3)\sqrt{v_1^2+v_2^2-k^2}}}{2 \sqrt{v_1^2+v_2^2-k^2}}, & v_1^2+v_2^2 >k^2,
            \end{dcases}
\end{align*}
and 
\begin{align*}
 &\int_0^\infty  K_0\big[(z_2^2+(x_3-y_3)^2)^{1/2}(v_1^2-k^2)^{1/2} \big] \cos(z_2v_2) \d z_2  \\
& =            \frac{\pi e^{(x_3-y_3)\sqrt{v_1^2+v_2^2-k^2}}}{2\sqrt{v_1^2+v_2^2-k^2}}, \quad  v_1^2+v_2^2 >k^2.
\end{align*}
Therefore, for $y_3>x_3$, we obtain
\begin{align*}
I_1 = \begin{dcases}
             \frac{i e^{-i(x_3-y_3)\sqrt{k^2-v_1^2-v_2^2} }}{2\sqrt{k^2-v_1^2-v_2^2}}, & v_1^2+v_2^2 <k^2, \\
            \frac{e^{(x_3-y_3)\sqrt{v_1^2+v_2^2-k^2}}}{2\sqrt{v_1^2+v_2^2-k^2}}, & v_1^2+v_2^2 >k^2.
\end{dcases}
\end{align*}
Similarly we have for $x_3+y_3<0$ that
\begin{align*}
I_2 &=  \int_{\R^2} \Phi(z_1,z_2,x_3+y_3)e^{-i(z_1v_1+z_2v_2)} \d z_1\d z_2 \\
 &=  \begin{dcases}
             \frac{ie^{-i(x_3+y_3)\sqrt{k^2-v_1^2-v_2^2} }}{2\sqrt{k^2-v_1^2-v_2^2}}, & v_1^2+v_2^2 <k^2, \\
            \frac{e^{(x_3+y_3)\sqrt{v_1^2+v_2^2-k^2}}}{2\sqrt{v_1^2+v_2^2-k^2}}, & v_1^2+v_2^2 >k^2.
            \end{dcases}
\end{align*}
Using $I_1$ and $I_2$ we derive for $x_3<0$ that
\begin{align*}
\widehat{w}(v_1,v_2,x_3) = 
\begin{dcases}
             \widehat{\varphi}(v_1,v_2)e^{-ix_3\sqrt{k^2-v_1^2-v_2^2}}, & v_1^2+v_2^2 <k^2, \\
          \widehat{\varphi}(v_1,v_2)e^{x_3\sqrt{v_1^2+v_2^2-k^2}}, & v_1^2+v_2^2 >k^2.
\end{dcases}
\end{align*}
Now by continuously extending $\widehat{w}(v_1,v_2,x_3)  = \widehat{\varphi}(v_1,v_2)$ when $v_1^2+v_2^2 =k^2$ we 
obtain~\eqref{eq:uhat}.
\end{proof}

\section{Experimental data: measurement and preprocessing}

\label{sect:imag} 

In this section we first give a brief description of our experimental setup.
Second, we describe the data preprocessing procedure which is one of the most
important ingredients of the paper.

\subsection{Data collection}

\label{sect:5.1}

\begin{figure}[h!]
\centering
\subfloat[A photograph of our experimental setup]{\includegraphics[width=0.4\textwidth]{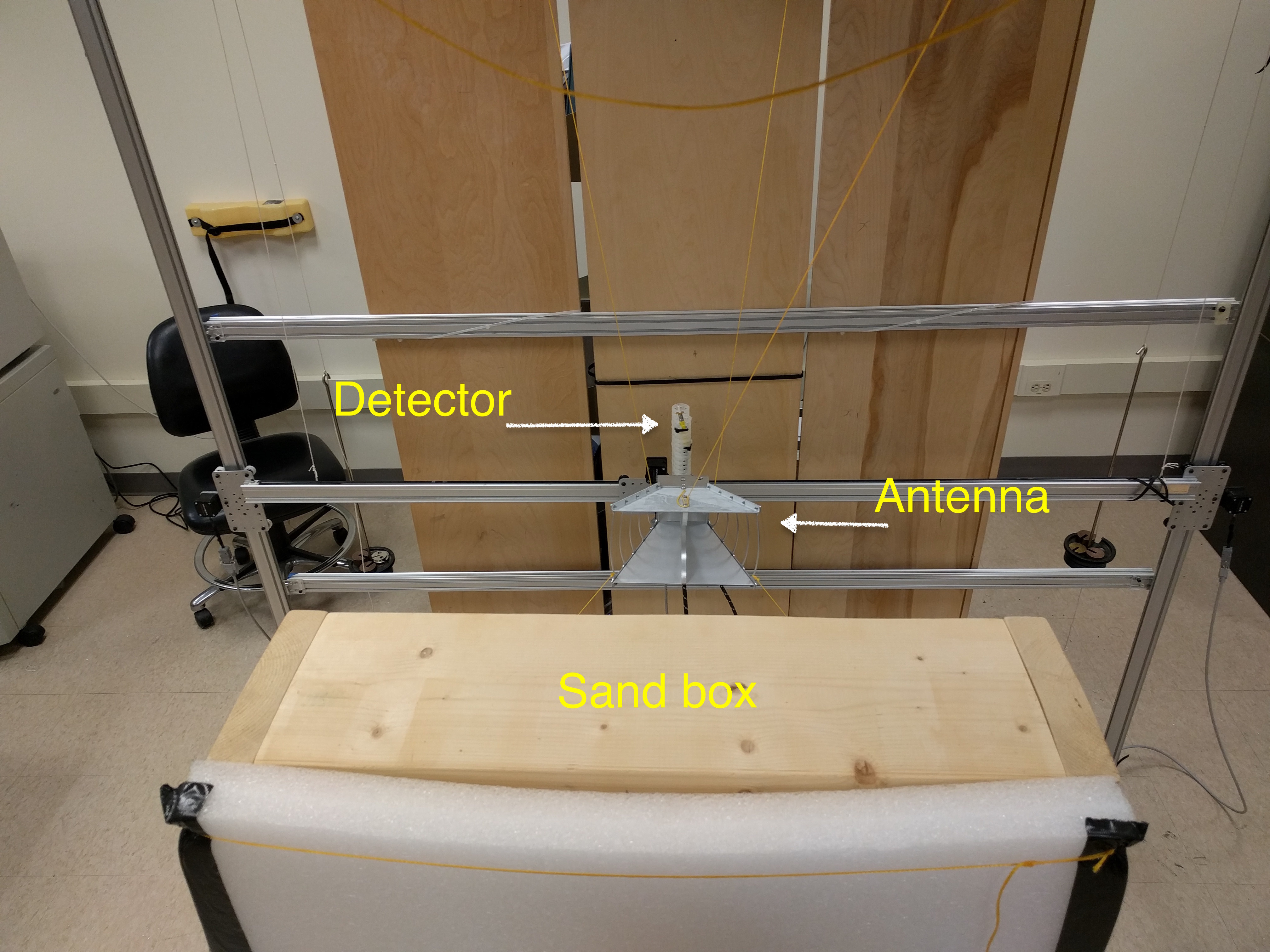}}\hspace{0.5cm} 
\subfloat[The schematic diagram of our experimental setup]{\includegraphics[width=0.5\textwidth]{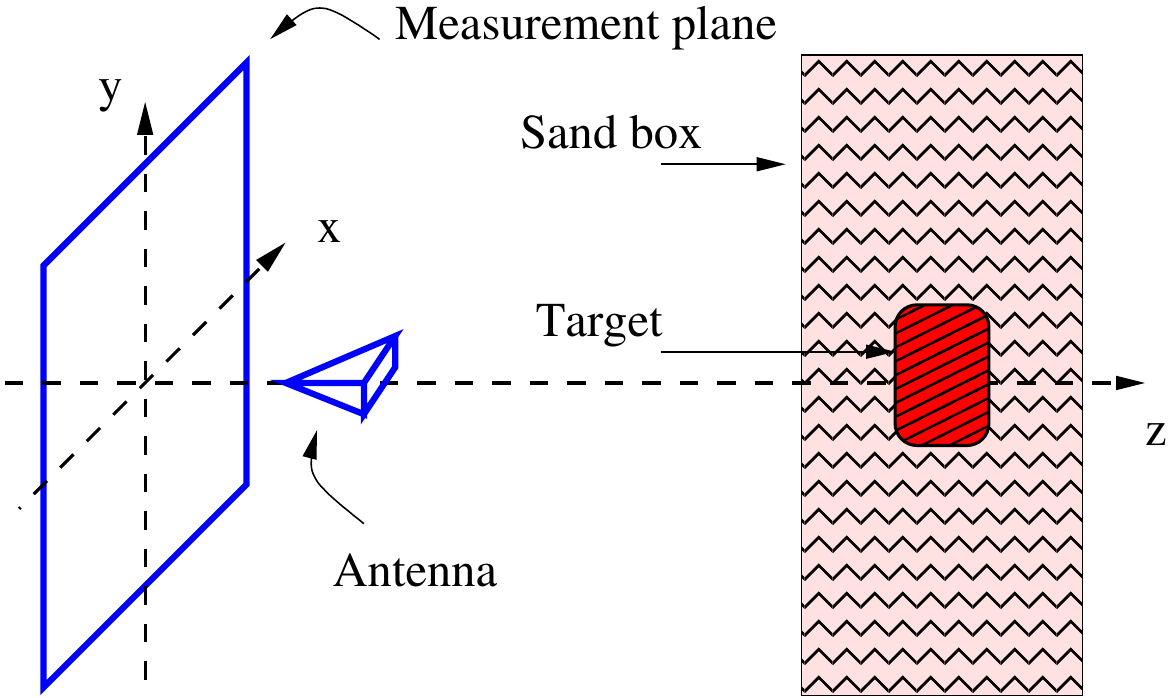}}
\caption{The experimental setup.}
\label{fi:1}
\end{figure}


The sandbox was filled with dry sand with no moisture in it. The
part of the plane on which our data are collected is a 1 meter by 1 meter
square. The $z-$axis is perpendicular to the measurement plane, while the $x-
$axis and the $y-$axis are the horizontal and the vertical axis,
respectively. The direction from the measurement plane to the target is the
positive direction of the $z-$axis. 

The broadband horn antenna and the collecting probe (detector) are both
connected to a 2-port Rohde $\&$ Schwarz vector network analyzer (VNA). The
distance between the horn antenna and the measurement plane was about 20 cm and
the distance between the measurement plane and the sandbox was about 75 cm. Due to
some technical reasons we placed the antenna in front of the measurement
plane. The VNA sends single frequency signals at frequencies ranging from 1
GHz to 10 GHz via the horn antenna. 
We collected the scattered wave at 300 frequency points uniformly
distributed over the range from 1 GHz to 10 GHz. Our data analysis
has shown that we can work only in a narrow interval of frequencies centered
at 2.6 GHz, 3.01 GHz or 3.1 GHz (section \ref{sect:centralFreq}). Since 2.6 GHz corresponds to the
wavelength of 11.5 cm and the burial depth of each target was
about 15 cm, then the distance between the horn antenna and each target was
about 71 cm $\approx$ 6.17 wavelengths. This distance can be considered as a sufficiently large
distance in terms of wavelengths. Thus, it justifies our modeling of the source
as a plane wave.

We have measured the backscatter data for a single direction of the incident
plane wave by repeating the experiment for different positions of the probe
on the measurement plane. More precisely, the probe was uniformly moved over
the scanning area with the step size 2 cm, that means a total of 2550 data
points for each of 300 frequencies. For the convenience in our numerical
implementation we work with dimensionless variables $\mathbf{x}^{\prime }=%
\mathbf{x}/(10\,\mathrm{cm})$ and keeping the same notation for brevity.
Thus, from now on, for example, $0.5$ of length means $5$ cm. In
principle, we can write equation (\ref{eqn Hel}) in terms of parameters with
dimensions via using $\mathbf{x}$ in cm and replacing $k$
 with the usual wavenumber $\tilde{k} =2\pi /\lambda$, 
 where $\lambda $ is the wavelength measured in centimeters. Next,
replacing $\mathbf{x}$ with $\mathbf{x}^{\prime }$ we obtain from (\ref{eqn Hel})
that the dimensionless wavenumber $k= \tilde{k} \cdot 10cm$.

Let $\mathbf{E}=(E_{x},E_{y},E_{z})$ be the electric field. The component $%
E_{y}$, which is the voltage, was incident upon the medium and the
backscatter signal of the same component was measured. Therefore, we denote $%
u(\mathbf{x},k):=E_{y}(\mathbf{x},k)$, where the function $u(\mathbf{x},k)$
is the above discussed solution of the problem~\eqref{eqn Hel}--\eqref{eqn outgoing}. 
Thus, the incident wave field is modeled as the plane
wave $E_{y,\mathrm{inc}}=e^{ikz}$. 
%

\subsection{Data preprocessing}

\label{sect:dataProcess} 

One of the challenges in working with the experimentally measured data is
that the data are perturbed by a significant amount of noise. We refer to
Figure~\ref{fi:2} for an example of our measured data at frequency 2.6 GHz.
The main reasons for the significant noise in the measurements are as follows:

\vspace{-0.0cm} 
\begin{figure}[h!]
\centering
\subfloat[Real part of measured data]{\includegraphics[width=0.4\textwidth]{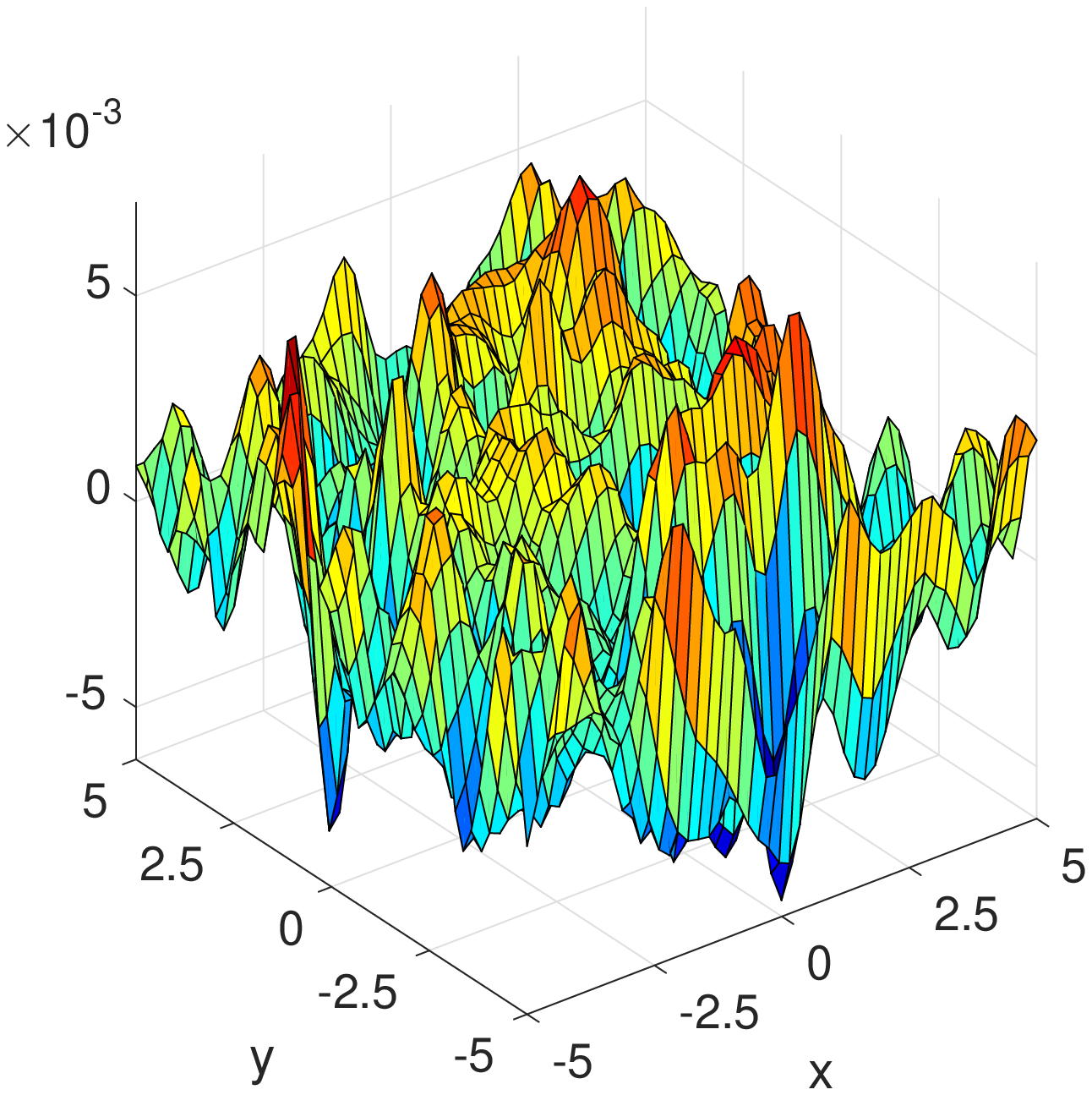}} \hspace{0.5cm} 
\subfloat[Imaginary part of measured data]{\includegraphics[width=0.4\textwidth]{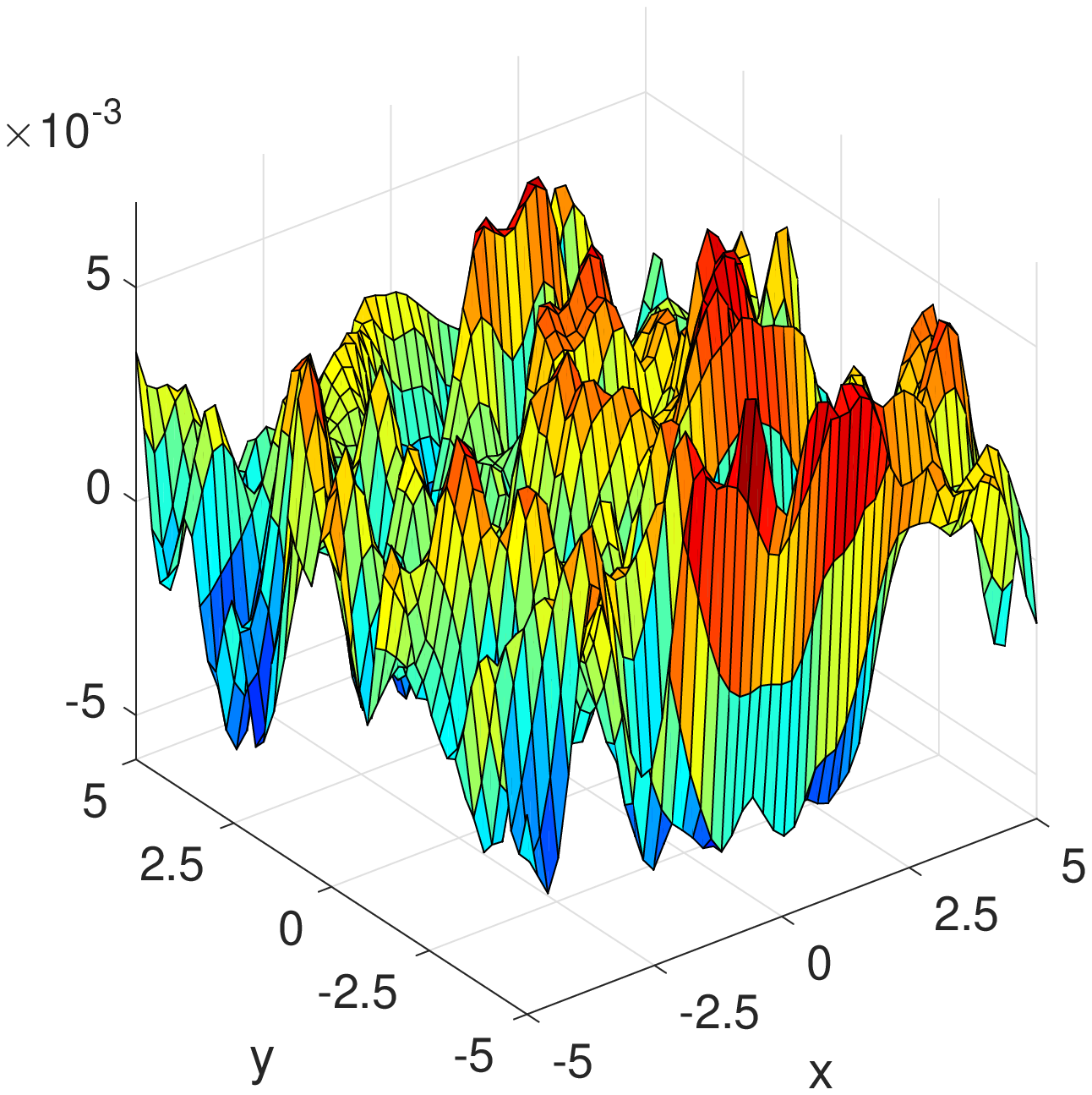}} 
\vspace{-0.4cm} 
\subfloat[Top view of (a)]{\includegraphics[width=0.35\textwidth]{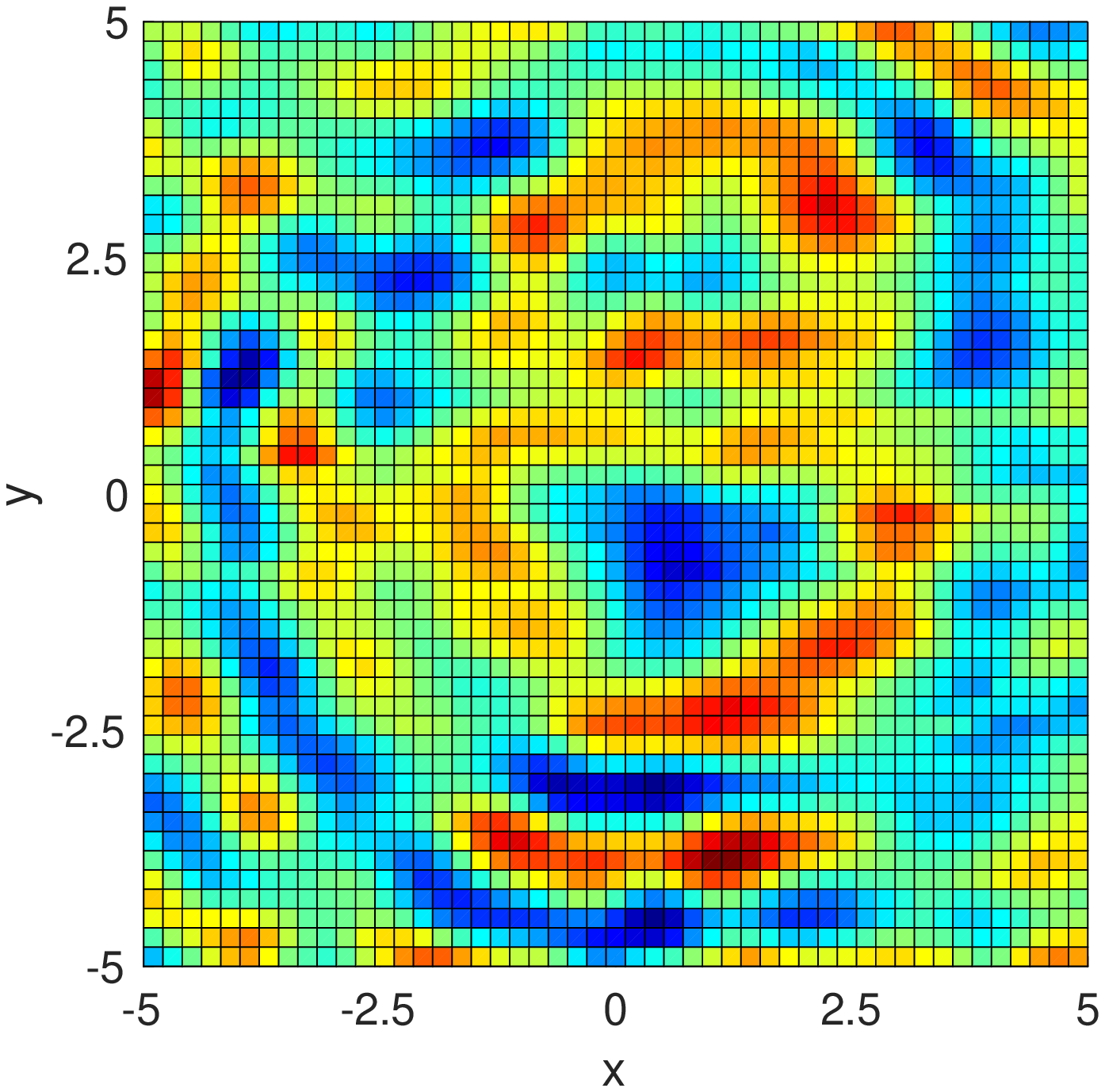}} \hspace{0.9cm}
\subfloat[Top view of (b)]{\includegraphics[width=0.35\textwidth]{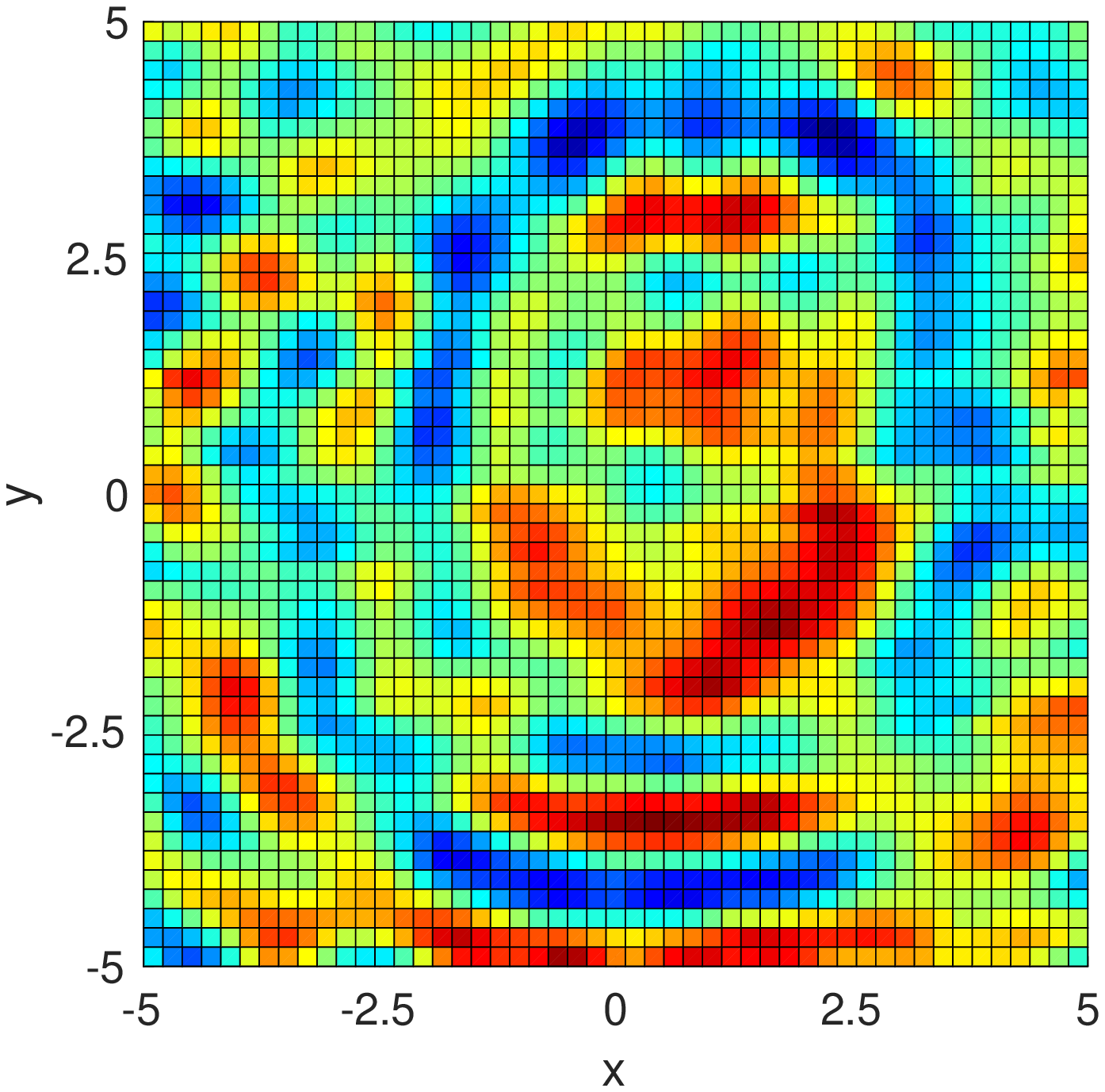}}
\caption{The yellow pine data on the measurement plane ($xy$-plane) at 2.6
GHz $(k=5.5)$.}
\label{fi:2}
\end{figure}

\begin{enumerate}
\item Due to some technical difficulties in the experimental setup, the
source antenna is placed between the measurement plane and the sandbox.
Hence, the waves reflected by the sandbox are scattered by the antenna
before they propagate to the probe. This causes some additional parasitic
signals.

\item Our experiments were conducted in a regular room, see Figure~%
\ref{fi:1}. Keeping in mind that in our desired application explosive-like
targets can be located in a cluttered environment, we intentionally did not
arrange sort of idealized conditions for our experiments, that means we did not
remove the regular furniture from this room and did not arrange a special
anechoic chamber. Therefore, our measured signals were reflected not only from the
sandbox and targets buried in it but also from other objects in the room,
such as the metallic bars placed behind the measurement plane and the room
furniture. Furthermore, these signals are affected by WiFi signals of the
building.

\item The incident waves emitted by the antenna propagate not only towards
the target but also partly backwards to the probe. 

\item The backscatter signal is unstable, due to the instability of the
emitted signal.
\end{enumerate}

Therefore, to apply our inversion method to the measured data, we need to
develop a heuristic data preprocessing procedure. The goal of this procedure
is twofold: (1) We somewhat distill the signals reflected by our buried
targets from signals reflected by the sandbox and other unwanted objects and
(2). We reduce the computational domain and the amount of noise in the
signal. Our data preprocessing procedure consists of three main steps:

\begin{enumerate}
\item[Step 1.] Given the data on the measurement plane for the case when a
target is buried in the sandbox, we subtract from it the reference data. 
The reference data are the ones, which were measured for the
case when the sandbox contains no buried targets. This subtraction helps us to sort
of extract the signals of the buried targets from the total signal and
also to reduce the noise.

\item[Step 2.] Back-propagating the data obtained after Step 1 using the
 data propagation procedure in Section~\ref{sect:SAR}. This data propagation aims to
 ``move" the data closer to the target. As a result, we
obtain reasonable estimates for the location of the buried targets,
particularly in the $(x,y)-$plane. In addition, this step helps us reduce
the computational domain $\Omega $.

\item[Step 3.] Determining an interval of frequencies on which the
propagated data are stable. We observe that this is a quite narrow interval.

\end{enumerate}

Next, the preprocessed data are used as the input for the algorithm
described in Section~\ref{sect:Algorithm1}. We point out once again
that it is our belief that even though the data preprocessing procedure is a
heuristic one (so as this is the common case in engineering), it is
justified in the end by the accuracy of our results, see Table 2 and Figures~\ref{fi:6} and~\ref{fi:7}.

\subsubsection{Subtraction and back-propagation}

\label{sec 4.2.1}

We recall that our inversion method is developed for the case when the
coefficient $c(\mathbf{x}) $ of equation (\ref{eqn Hel}) 
is a smooth function for $\mathbf{x}\in \mathbb{R}^{3}$. On the other
hand, since our targets are buried in a sandbox, we actually have air/sand
interface, on which $c(\mathbf{x}) $ has a discontinuity.
Clearly the mathematical model with the air/sand interface would be a more
complicated then the one with a smooth function $c(\mathbf{x}) .$
In addition, it is yet unclear how to modify our method for this case.
Hence, to avoid considerations of this model, we conduct an important step
of our data preprocessing procedure. Namely, we subtract the reference signal  
from the total signal, see Step 1 in Section~\ref{sect:dataProcess}  for the
definition of the reference signal. First, we  measured the data for the
sandbox without any buried objects in it. For each buried target and for
each location of the probe, we subtract the latter data from the data
measured when this target was present. The data obtained after this
subtraction are then treated  as the backscatter data from the targets being
placed in a medium with a smoothly varying function $c(\mathbf{x}).$

From now on, the term ``data" is applied only to the data
obtained by the above subtraction. The next step is to back-propagate the
data. We call this the data propagation procedure, see Section 4. Given the
data on the measurement plane, the goal of the data propagation is to
approximately  ``move" the data to a plane which is closer to
the unknown target than the measurement plane. The back-propagation first
allows us to reduce the computational domain for our algorithm. Second, it
makes the backscatter data more focused. Finally, we have also observed that
the back-propagation helps us to significantly reduce the noise in the
signals. These last two observations, particularly the focusing of the
propagated data, are crucial for estimating of the location of the targets
and further steps of the data preprocessing. 

As mentioned in Section~\ref{sect:SAR}, another name of this
back-propagation technique is \textquotedblleft the angular spectrum
representation" \cite[Chapter 2]{Novot2012}. The authors have
successfully used this technique in~\cite{Kliba2016,Exp1,Exp2}. We now specify
considerations of Section 4.1 for our experimental setup.

Let $a$ and $b$ be two numbers such that $b<a<0$. Let the domain of interest $\Omega \subset \left\{ z>a\right\} .$
Let $g(\mathbf{x},k)$ be the data defined on the measurement rectangle 
\begin{equation*}
P_{m}=\{\mathbf{x}:-5<x,y<5,z=b\}
\end{equation*}%
of the plane $\left\{ z=b\right\} $ and let  $f(\mathbf{x},k)$ be its
approximation on the propagated rectangle 
\begin{equation*}
P_{p}=\{\mathbf{x}:-5<x,y<5,z=a\}
\end{equation*}%
of the plane $\left\{ z=a\right\} $. Extend both functions $g$ and $f$ by
zero outside of the rectangles $P_{m}$ and $P_{p}$, respectively. For the
convenience of the presentation, we call below sets $P_{m}$ and $P_{p}$
\textquotedblleft measurement plane" and ``propagated plane"
respectively. Note that since  $\Omega \subset \left\{ z>a\right\} ,$
then $P_{p}$ is closer to the target than $P_{m}$.

\begin{remark}
\label{convention}
We note that the data propagation formula~\eqref{4.12} is derived under a convention that $e^{ikz}$
is a plane wave propagating to the positive direction. It seems to us that this convention is 
popular in the community we know.  We present our theory for this convention for 
the convenience of the readers. However, it is known there also exists a convention that $e^{-ikz}$ is 
a plane wave propagating to the positive direction. Our theory in this paper of course can be adapted
easily to this convention. 
The formula~\eqref{4.12} has been successfully implemented for simulated data, see~\cite{Kliba2016} or Figure~\ref{fi:5}(b). 
However, we have noticed that for our experimentally measured data here, to obtain positive results we must use the data
propagation formula for the convention that the plane wave $e^{-ikz}$ propagates
to the positive direction. We believe that this is the convention for the machine
in our experimental setup.  
\end{remark}
From the above remark, we modify the data propagation formula for our experimental data as follows:
\begin{equation}
f(x,y,a,k)=\frac{1}{2\pi }\int\limits_{v_{1}^{2}+v_{2}^{2}<k^{2}}\widehat{g%
}(v_{1},v_{2},k)e^{i\left( xv_{1}+yv_{2}\right) }e^{-i\sqrt{k^{2}-v_{1}^{2}-v_{2}^{2}} \left(
b-a\right) }\d v_{1}\d v_{2},  \label{eq:propdata}
\end{equation}%
where 
\begin{equation*}
\widehat{g}(v_{1},v_{2},k)=\frac{1}{2\pi }\int\limits_{\mathbb{R}%
^{2}}g(x,y,b,k)e^{-i(xv_{1}+yv_{2})}\d x \d y.
\end{equation*}

\begin{figure}[h]
\centering
\subfloat[Before data
propagation]{\includegraphics[width=0.4\textwidth]{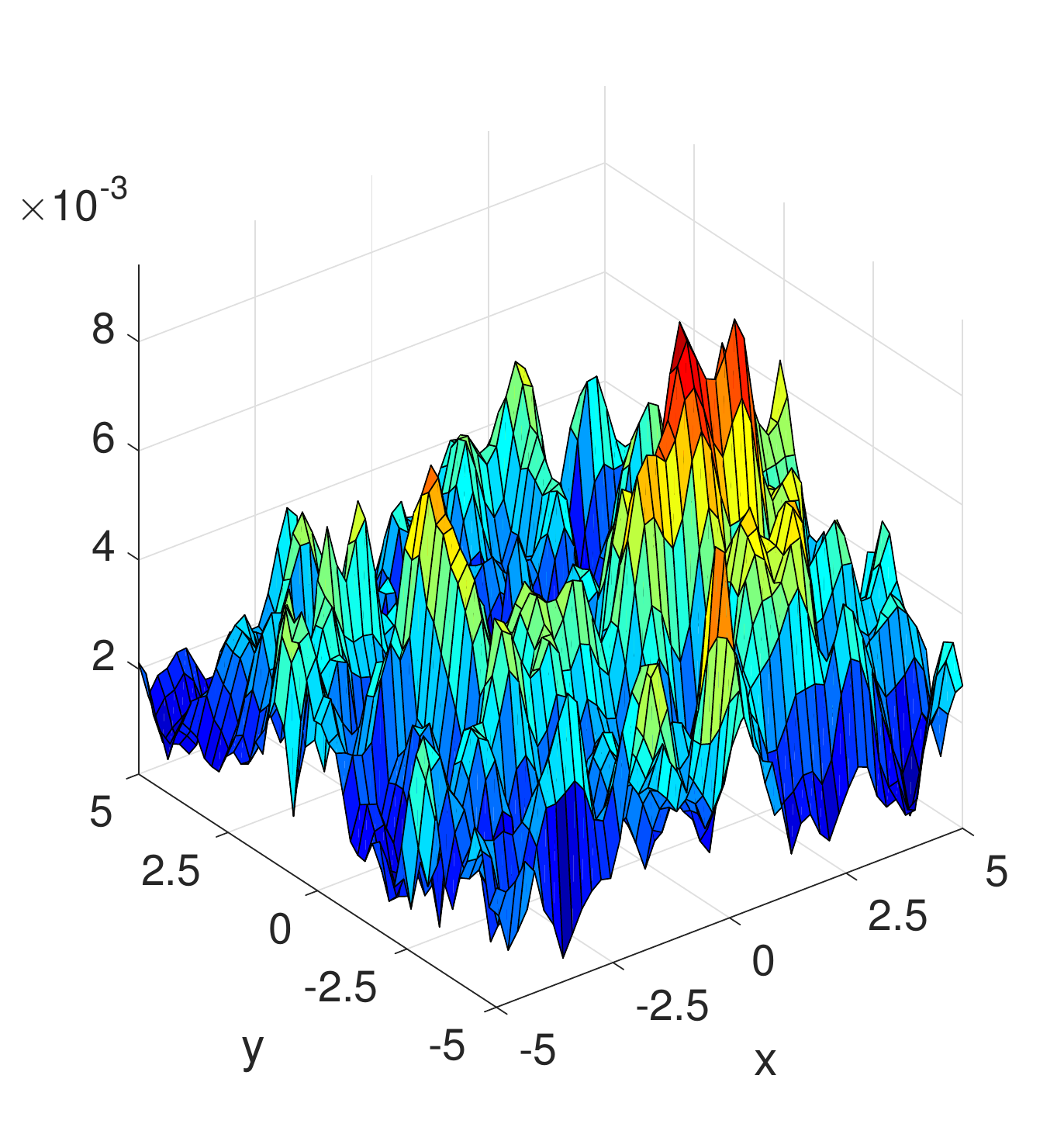}} 
\hspace{0.5cm} 
\subfloat[After data
propagation]{\includegraphics[width=0.4\textwidth]{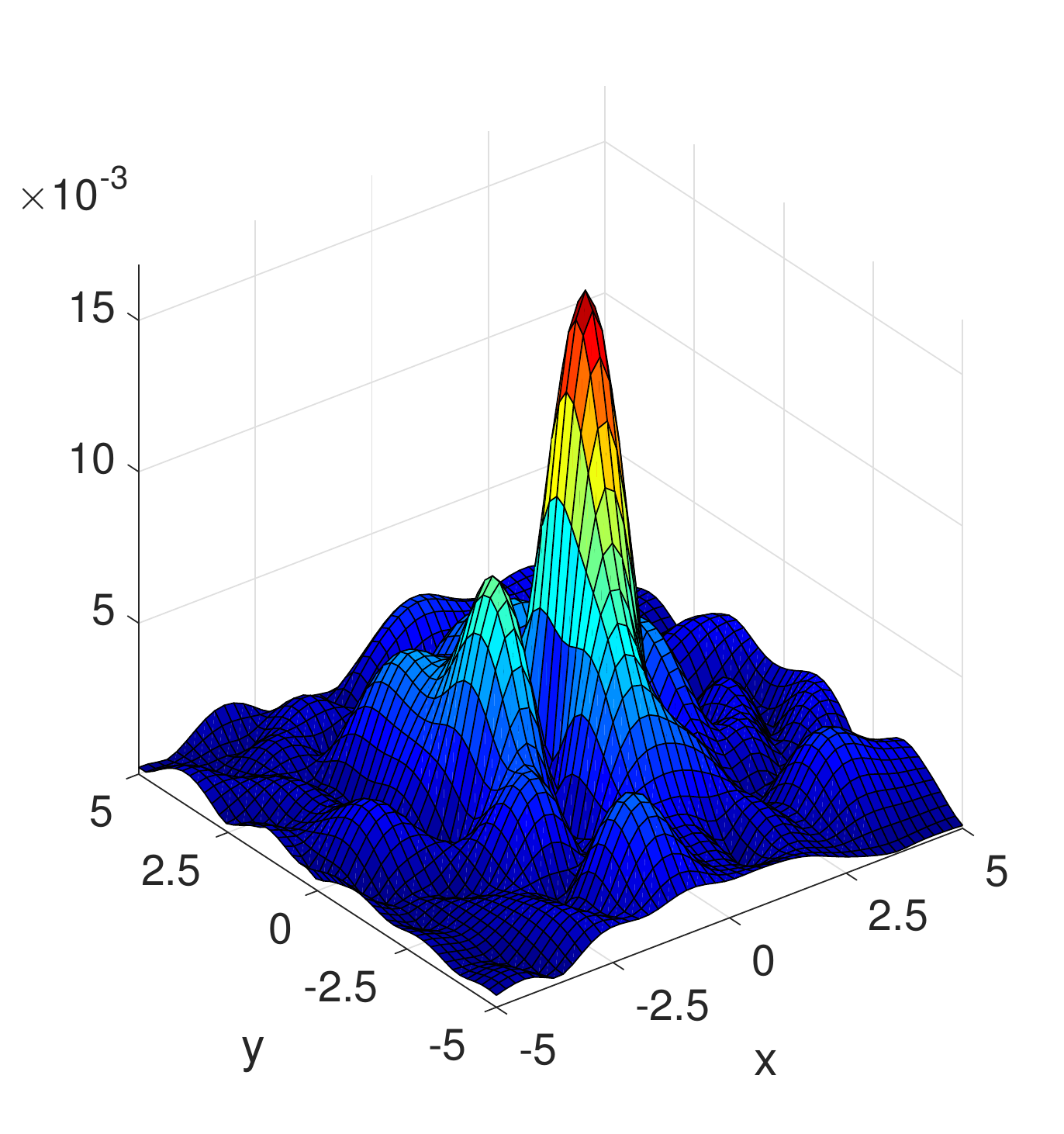}}
\caption{We respectively display in (a) and (b) the modulus of the (yellow
pine) data before and after back-propagation for $k=5.5$. }
\label{fi:3}
\end{figure}

In Figure~\ref{fi:3}(b) we present the absolute value of the experimental
data for one of our targets on a rectangle $P_{p}$ of the propagated plane
$P_{p}=\{\mathbf{x}:-5<x,y<5,z=-0.8\}$,
computed by the formula~\eqref{eq:propdata} for the wavenumber $k=5.5$. 
The corresponding rectangle on the measurement plane was
$P_{m}=\{\mathbf{x}:-5<x,y<5,z=-8.78\}$. Comparing Figure~\ref{fi:3}(a) with 
Figure~\ref{fi:3}(b), one can easily see
that the focusing of the scattered field is considerably improved after data
propagation. A similar behavior was observed for all other buried
targets we have worked with. From now on we will be interested only in the
propagated data instead of the original one. 

\subsubsection{Choosing an appropriate frequency interval and data truncation}

\label{sect:centralFreq}
\begin{figure}[h!]
\centering
\subfloat[]{\includegraphics[width=0.5\textwidth]{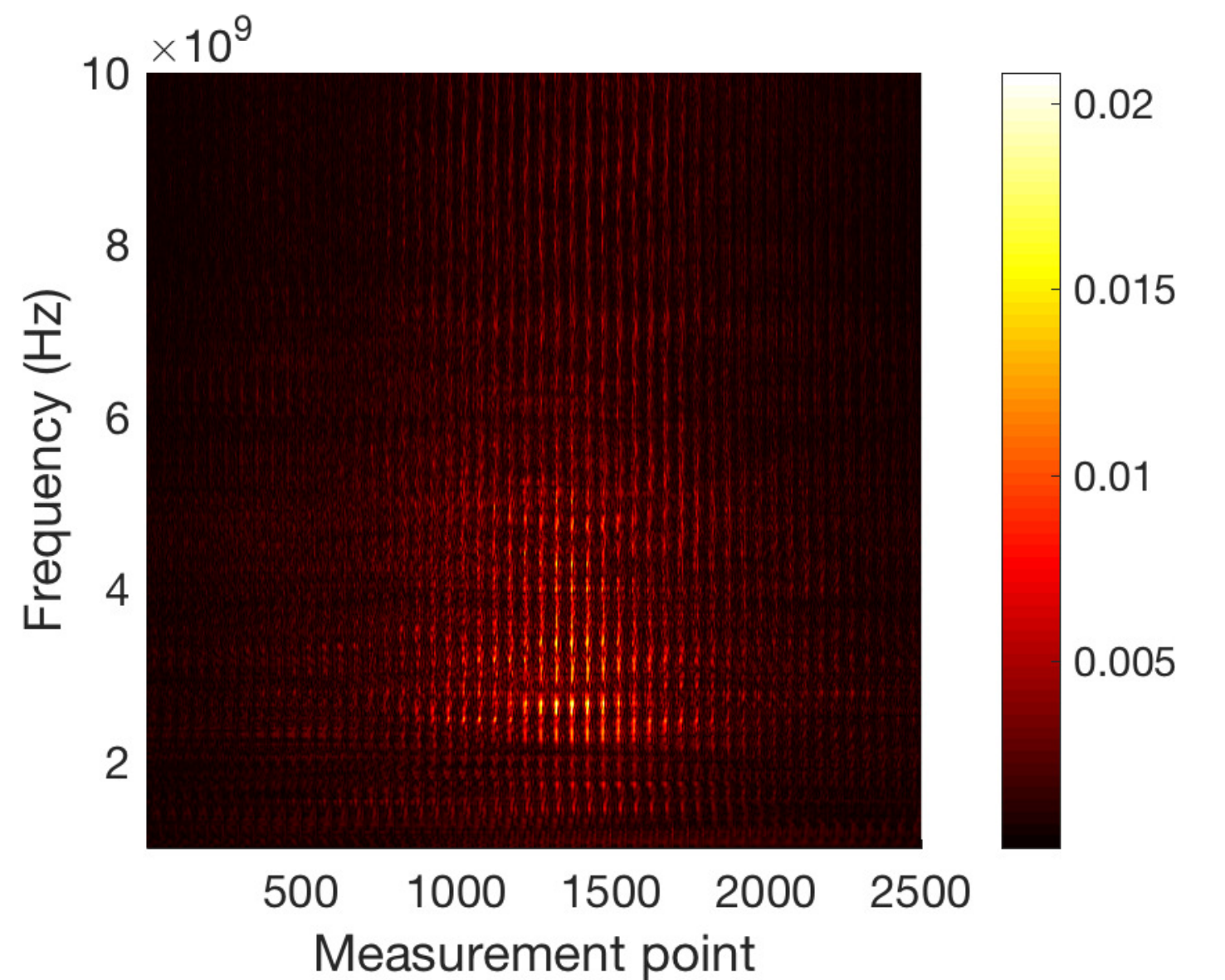}} \hspace{0.5cm} 
\subfloat[]{\includegraphics[width=0.43\textwidth]{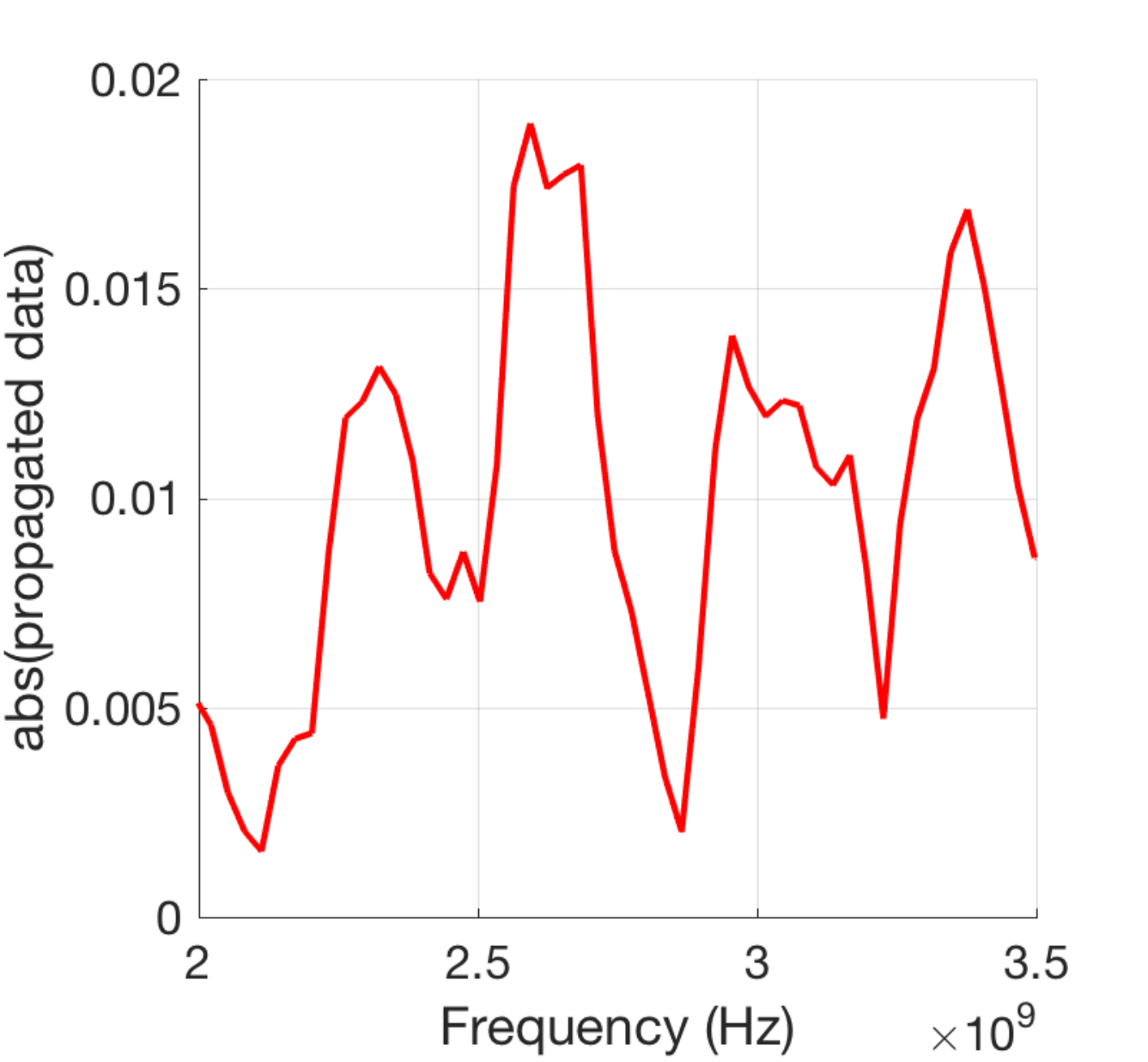}}\\
\caption{The absolute value of the experimental data after propagation for
all frequencies ranging from 1 GHz to 10 GHz. (a) Top view; (b) Restriction on the interval of frequencies $[2, 3.5]$ GHz for 
 the 1373th measurement point.}
\label{fi:4}
\end{figure}
We first describe how to choose an interval of frequencies where the
propagated data are usable for our inversion. Figure~\ref{fi:4}(a) presents
the absolute value of the propagated data for all 300 frequencies from 1 GHz
to 10 GHz. The propagated data are written as a matrix of $300\times 2550$.
We can roughly see from Figure~\ref{fi:4}(a) that these signals are in general stronger and more
focused for frequencies from 2 GHz to 3.5 GHz. The signals for other
frequencies are weak and they are dominated by the noise. Therefore, they
are not usable for the inversion. Now if we observe more carefully in Figure~\ref{fi:4}(b), 
the stronger signals are not stable on the entire interval of frequencies $%
[2,3.5]$ (GHz). Indeed, Figure~\ref{fi:4}(b) shows the behavior of the modulus of the propagated
data at the 1373th measurement point. It can be seen from this figure that
that the maximum of the modulus of the propagated data can change quite
rapidly over some very small intervals of frequencies.


We are interested in  an interval of frequencies where the data satisfy
the following two criteria: (1) In this interval, the maximal value of the
modulus of the propagated data does not change too rapidly, and (2) For
different frequencies within this interval, these maxima are attained at the
same coordinates (up to small error) on the propagated plane.

Following the two criteria one can easily determine such intervals of
frequencies. We choose the longest interval of frequencies with data
satisfying those two criteria. For the case in Figure~\ref{fi:4}(b), 
$[2.53,2.71]$ (GHz) is our stable interval of frequencies. Now
if the stable interval of frequencies has an even number of frequencies, we define that the optimal frequency
is the smaller number of the two middle numbers of the interval,
otherwise the optimal frequency is the median of the interval. The optimal frequency
for $[2.53,2.71]$ (GHz) is $2.62$ GHz. 
 The corresponding interval of wavenumbers is $[\underline{k},\overline{k}]=\left[ 5.31,5.69
\right] $. We observed that the lower and upper bounds of such intervals
differ slightly for the data from different targets considered in this
paper. From now on, by using $[\underline{k},\overline{k}]$ we understand
that this interval of wavenumbers corresponds to the interval of frequencies
which was chosen above. Thus, we have stable propagated data on this
interval.

\begin{figure}[!h]
\centering
\subfloat[Propagated data after
truncation]{\includegraphics[width=0.4\textwidth]{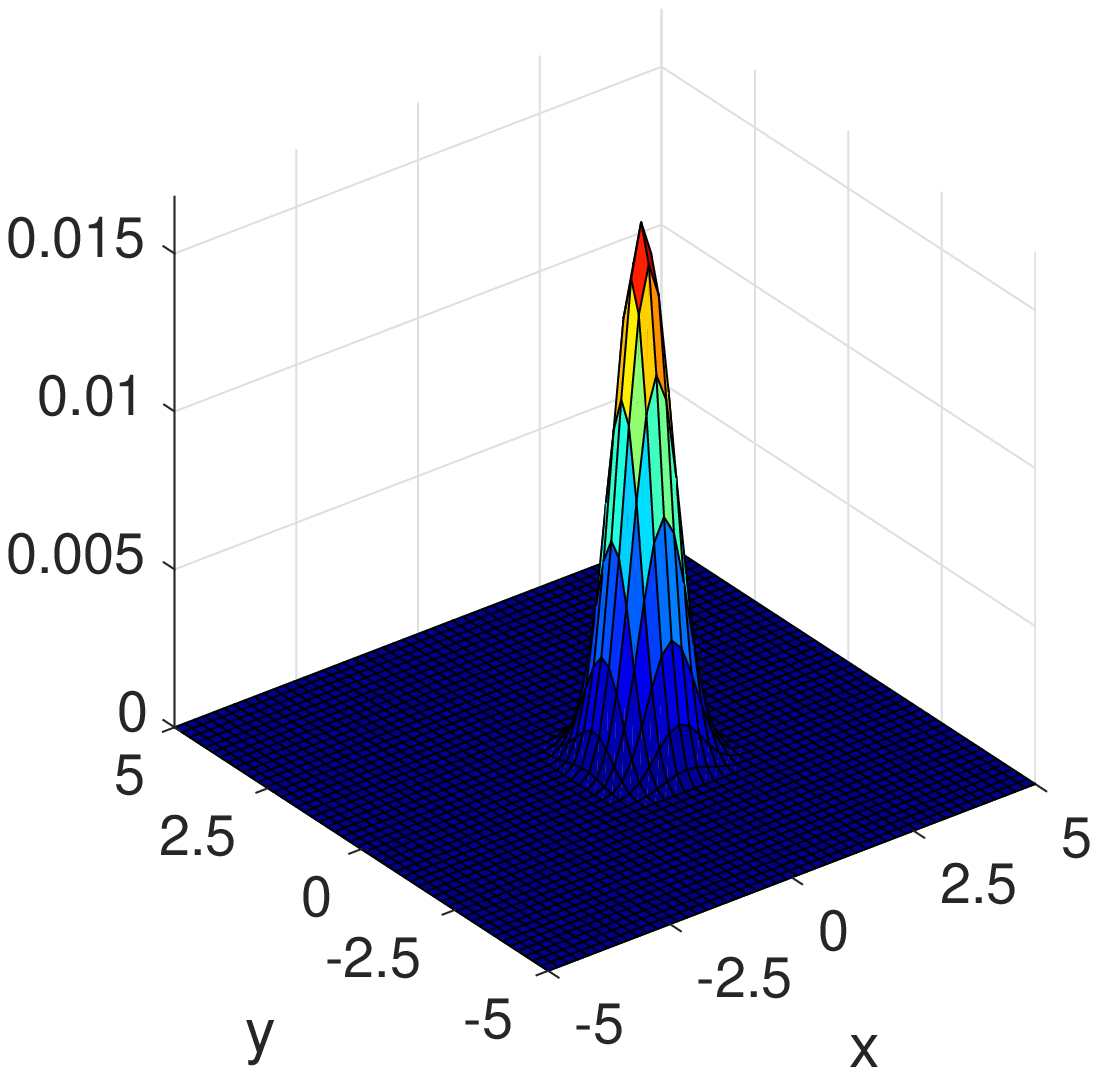}} 
\hspace{0.5cm} 
\subfloat[Propagated data in a
simulation]{\includegraphics[width=0.4\textwidth]{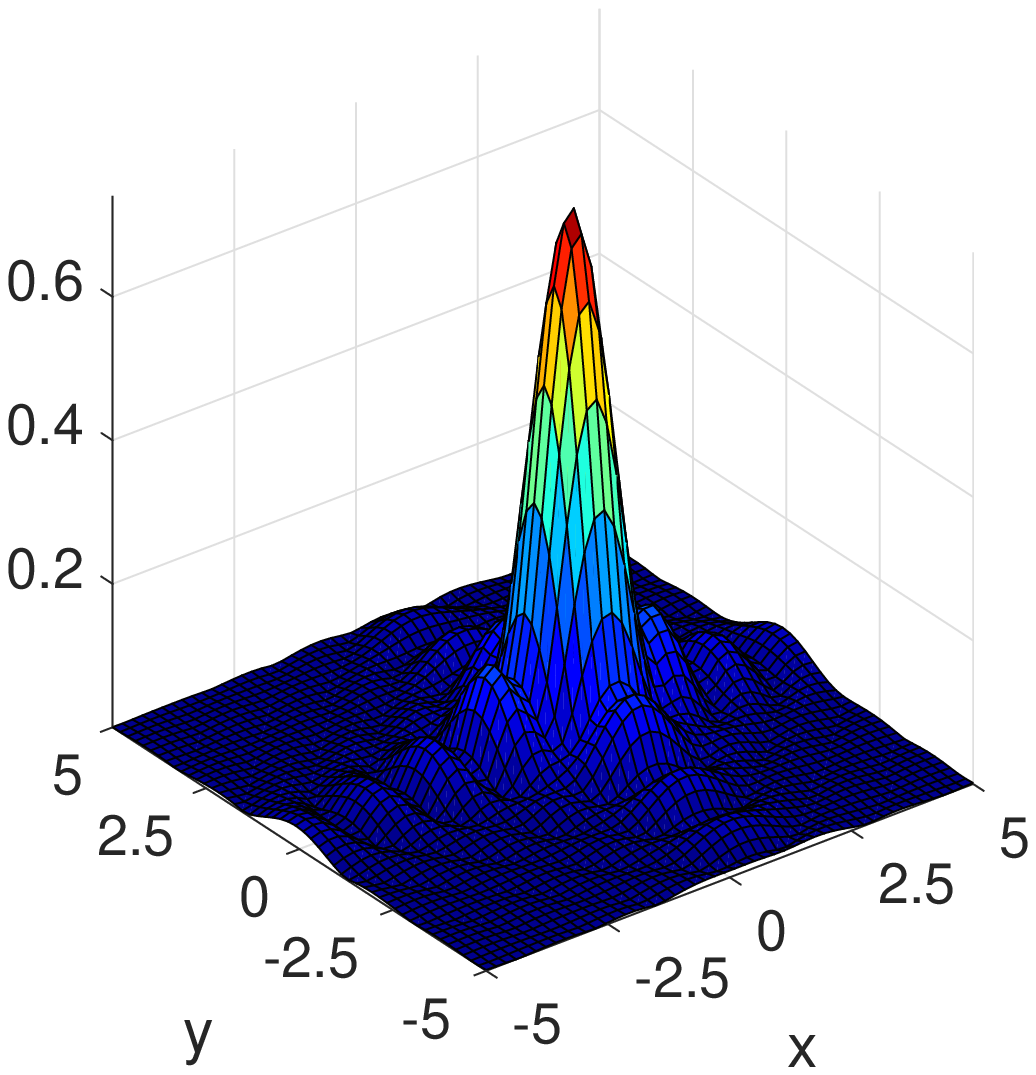}}
\caption{The absolute value of the propagated data after truncation and the
propagated data in simulation for $k=5.5$. The peaks in the both figures 
provide a good estimation of the $xy$-location of the targets. Truncation was not
performed for simulated data of (b).}
\label{fi:5}
\end{figure}

We now describe the last step of the data preprocessing procedure called
\textquotedblleft data truncation". It can be seen from the propagated data
that although its magnitude always attains the global maximum at the $xy-$%
location of the target for all $k\in \lbrack \underline{k},\overline{k}]$,
there are oscillations on other parts of the propagated plane, see Figure~%
\ref{fi:3}(b). These oscillations seem to be random, i.e. they are all
different for different frequencies. Hence, we conclude that the propagated
data are dominated by noise outside a certain neighborhood of their global
maxima. We believe that this noise can be reduced by \textquotedblleft
cleaning" these random oscillations. This means that we truncate the
propagated data as follows. For each $k\in \lbrack \underline{k},\overline{k}%
]$ we replace the function $f(x,y,a,k)$ with the function $\widetilde{f}%
(x,y,a,k),$ where 
\begin{equation}
\widetilde{f}(x,y,a,k)=%
\begin{cases}
f(x,y,a,k), & \text{ if }|f(x,y,a,k)|\geq 0.8\max_{x,y}|f(x,y,a,k)|, \\ 
0, & \text{ otherwise.}%
\end{cases}
\label{eq:truncation}
\end{equation}%
Next, the function $\widetilde{f}(x,y,a,k)$ is smoothed by a Gaussian
filter. Numerically the smoothing function \texttt{smooth3} in MATLAB with
the Gaussian option is used. This is done for all values of $k\in \lbrack 
\underline{k},\overline{k}]$. We display in Figure~\ref{fi:5} the propagated
data after the truncation for $k=5.5$. It can be seen that the behavior of
the absolute value of these processed data looks more similar to that of
simulated data. The results are similar for other wavenumbers $k$ in $[%
\underline{k},\overline{k}]$.

\section{Numerical implementation and reconstruction results}

\label{sect:result}

We describe in this section some details of the numerical implementation and
present the reconstruction results for our experimental data using the
globally convergent algorithm. We test here six data sets. Each data set
corresponds to a single scattering object, numbered from 1 to 6, see Table~%
\ref{tab:table1}. \vspace{-0.0cm} 
\begin{table}[tbph]
\caption{Scattering objects}
\label{tab:table1}\centering
\begin{tabular}{|c|c|}
\hline
Object ID & Scattering object \\ \hline
1 & A piece of bamboo \\ \hline
2 & A geode \\ \hline
3 & A piece of rock \\ \hline
4 & A piece of sycamore \\ \hline
5 & A piece of wet wood \\ \hline
6 & A piece of yellow pine \\ \hline
\end{tabular}%
\end{table}

\subsection{Target location estimation}

\label{sect:6.1}

\begin{figure}[h!]
\centering
\subfloat[Yellow pine
data]{\includegraphics[width=0.4\textwidth]{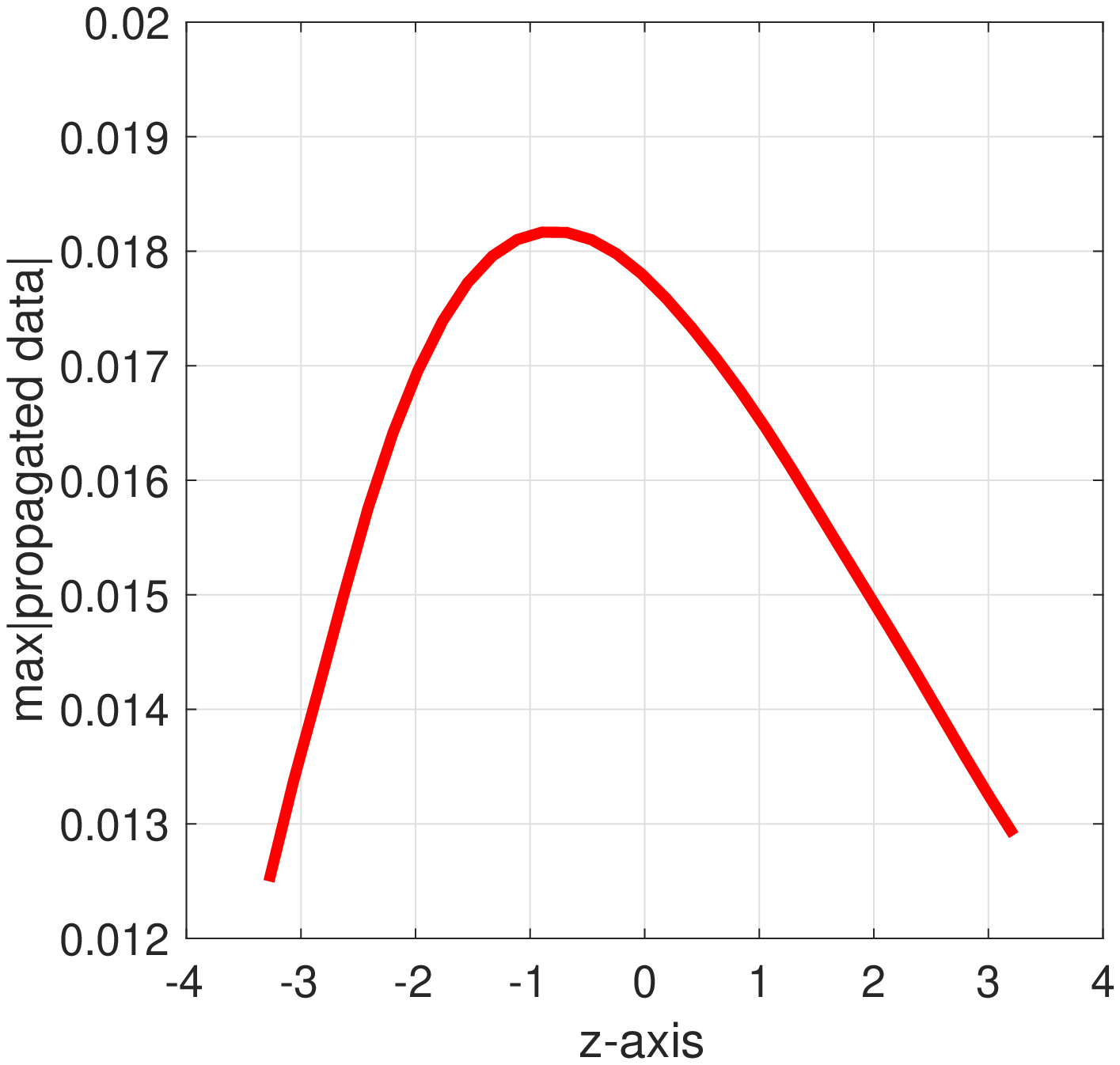}} \hspace{0.5cm}
\subfloat[Geode data]{\includegraphics[width=0.44\textwidth]{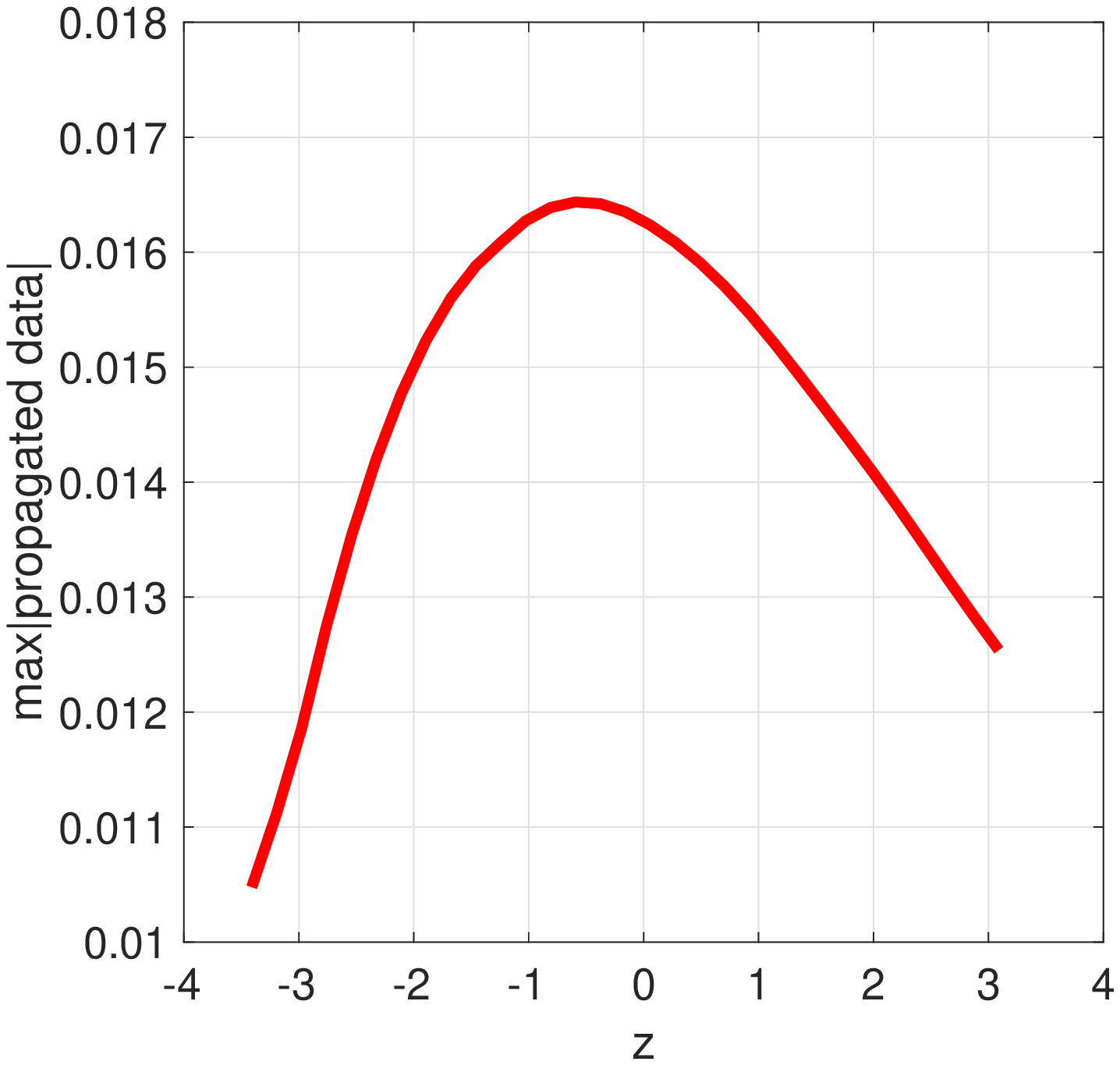}}
\caption{We respectively display in (a) (yellow pine data) and (b) (geode
data) the maxima of modulus of the propagated data with respect to different
propagated planes.}
\label{fi:5a}
\end{figure}

First, as in Section~\ref{sect:centralFreq}, for each target we can choose the frequency 
for which the signals of the propagated data are the strongest ones, see  Figure~\ref{fi:4}. 
In particular, in numerical studies, which correspond to Figure~\ref{fi:5a}, we use the
frequency of 2.6 GHz, i.e. $k=5.46$. For example, Figure~\ref{fi:5a}
presents the curves of maximal values with respect to $x,y$ of the modulus
of the propagated data for different propagated planes $\{z=a\}$, where $a$
varies between $-3.5$ and $3.5$. Since the thickness of the sandbox
in the $z-$direction was 50 cm (5 in dimensionless coordinates),
then this range is long enough to include the buried targets. We have
observed that these curves attain their maximal values at $z=z^{\ast }$
where $z^{\ast }<0$ and the distance between $z^{\ast }$
 and the front face of the target varies between $0.4$
and $1$ for all six data sets. Hence, the number $z^{\ast }$
provides us with a rough estimate of the burial depth of the front
face of that target, or, equivalently, an estimate of its $z-$coordinate.
This helps us substantially  reduce our search domain for the targets
in the $z$-direction. We refer to Section~\ref{sect:details} for details about
our search domain for the targets.

As already mentioned in Section~\ref{sect:centralFreq} and before, 
the data propagation provides an estimation for the location of a target in
the $xy$-plane. More precisly, we
define $\Omega _{T}$ as 
\begin{equation}
\Omega _{T}=\left\{ (x,y):|\widetilde{f}_{smth}(x,y,z^{\ast },\widetilde{k}%
)|>0.7\max |\widetilde{f}_{smth}(x,y,z^{\ast },\widetilde{k})|\right\}
\subset P_{p},  \label{OmegaT}
\end{equation}%
where $\widetilde{f}_{smth}(x,y,z^{\ast },\widetilde{k})$ is the propagated
data after the truncation procedure in the previous section and $\widetilde{k%
}$ is the wavenumber corresponding to the optimal frequency determined in
Section~\ref{sect:centralFreq}. Note that for each target $|\widetilde{f}%
_{smth}(x,y,z^{\ast },k)|$ has a positive peak whose location is the same
for all $k\in \lbrack \underline{k},\overline{k}]$, see Figure~\ref{fi:5}.
The truncation value 0.7 was chosen based on trial-and-error tests on
simulations. We have observed that $\Omega _{T}$ provides a reasonable
approximation for the $xy-$location of a target. The same truncation was
applied to all targets. Hence, it is not biased.

\subsection{Some computational details}

\label{sect:details} We consider the computational domain $\Omega $ as 
\begin{equation}
\Omega =(-2.5,2.5)\times (-2.5,2.5)\times (z^{\ast },5+z^{\ast }),
\label{6.1}
\end{equation}%
where the number $z^{\ast }<0$ is found as in the previous section.
The front face of the target is positioned at $\{z=0\}$ in our
numerical implementation. The data are propagated to the rectangle $P_{p}$
located on the plane $\{z=z^{\ast }\}$ and the unknown
targets are searched in the range $(z^{\ast },1)$ of the $z$-direction. 
The motivation of this search range is from similar
ranges that have been used in~\cite{Thanh2014} and our desired
applications for detection of buried explosives. Indeed, $z\in (z^{\ast },1)$ 
means that the linear size of the target in the $z$-direction does not exceed $1+z^{\ast }$ which is about between 5
cm and 15 cm in our cases. It is known that the latter range is a typical
linear size of antipersonnel mines and IEDs. 
We note that the larger range $(z^{\ast },5+z^{\ast })$ is
considered for the goal of completing the backscatter data in the next
section. Also in the $xy-$plane we restrict ourselves to the smaller area $%
(-2.5,2.5)^{2}$ than the original measurement one $(-5,5)^{2}$. The main
reason comes from the observation that the truncated data is zero outside of
the area $(-2.5,2.5)^{2}$, see Figure~\ref{fi:5}, and does not contribute to
the reconstruction process. This restriction helps us reduce our
computational domain.

We solve the Lippmann-Schwinger integral equation~\eqref{eqn LS} using a
spectral method developed in~\cite{Lechl2014}. The method relies
on a periodization technique of the integral equation, which enables the use
of the fast Fourier transform in the numerical implementation. The boundary
value problem~\eqref{30}--\eqref{32} is solved by a finite element
method. The implementation was done in FreeFem++~\cite{Hecht2012}, a
standard software designed with a focus on solving partial differential
equations using finite element methods, see~\emph{www.freefem.org} for more
information about FreeFem++.

Keeping in mind that we seek for the target in $\left\{ z\in (z^*,1)\right\} $ and using $\Omega _{T}$ in~\eqref{OmegaT}, we truncate the
coefficient $c_{n,i}$ as follows 
\begin{equation}
\widetilde{c}_{n,i}(\mathbf{x}):=%
\begin{cases}
\max \left( |c_{n,i}(\mathbf{x})|,1\right) , & \mathbf{x}\in \Omega
_{T}\times (z^{\ast },1), \\ 
1, & \text{elsewhere}.%
\end{cases}
\label{6.2}
\end{equation}%
After that we smooth the function (\ref{6.2}) by \texttt{smooth3} in MATLAB
and then the so obtained function $c_{n,i}\left( \mathbf{x}\right) $ is used
for solving the Lippmann-Schwinger equation in Algorithm~\ref{alg:globalconv}%
.

\subsection{Data completion}

\label{sect:completion}

We recall that only the backscatter signals were measured in our
experiments. This means that after the above data preprocessing procedure,
the function $g(\mathbf{x},k)$ was known only on the side $\Gamma =\{\mathbf{%
x}\in \partial \Omega :z=z^{\ast }\}$ of the domain $\Omega $ defined in (%
\ref{6.1}). As in~\cite{Thanh2014, Kliba2016}, we extend the scattered field
by zero on the other parts of $\partial \Omega$. We believe that this extension 
is reasonable because these parts of the boundary $\partial \Omega$ are relatively far from our search 
domain of the targets, and the scattered fields are supposed to be dominated by the 
incident field on those boundaries. In other words, in our
computation, we consider $g(\mathbf{x},k)$ on the entire boundary $\partial
\Omega $ as 
\begin{equation}
g(\mathbf{x},k):=%
\begin{cases}
g(\mathbf{x},k), & \mathbf{x}\in \Gamma , \\ 
e^{-ikz}, & \mathbf{x}\in \partial \Omega \setminus \Gamma .%
\end{cases}
\label{eq:completion}
\end{equation}
The use of $e^{-ikz}$ instead of $e^{ikz}$ follows from the convention for the
experimental data discussed in the end of section~\ref{sec 4.2.1}.
It is known that data completion methods are widely used for inverse
problems with incomplete data. Our data completion in~\eqref{eq:completion}
is a heuristic one relying on the successful experiences of our group when
working with globally convergent methods for experimental data, see~\cite%
{Thanh2014, Beili2012a}. Other data completion methods may also be applied.

\subsection{The stopping rules and choosing the final result}

\label{sect:stopping}

We recall the content of convergence theorem 6.1 of~\cite{Kliba2016} that
computed coefficients $c_{n,i}(\mathbf{x}) $ are located in a
sufficiently small neighborhood of the exact coefficient, provided that the
number of iterations is not too large. The theorem does not state that these
computed coefficients tend to the exact one, also, see theorem 2.9.4 in~\cite%
{Beili2012} and theorem 5.1 in~\cite{Beili2012a} for similar results.
Therefore, we address the stopping rules computationally.

These rules are about the choice of the final result for our Algorithm~\ref%
{alg:globalconv}. They have been introduced in~\cite{Exp1}. They rely on the
content of the convergence theorem of~\cite{Kliba2016} and trial-and-error
testing for simulated data. For the convenience of the readers, we present
these rules here.

We have two stopping rules: one for the inner iterations and the second one
for the outer iterations. Before stating those criteria precisely we need
some definitions. Denote by $e_{n,i}$ the relative error between the two
computed coefficients corresponding to two consecutive inner iterations of
the $n$-th outer iteration. In other words, 
\begin{equation*}
e_{n,i}=\frac{\Vert c_{n,i}-c_{n,i-1}\Vert _{L_{2}(\Omega )}}{\Vert
c_{n,i-1}\Vert _{L_{2}(\Omega )}},\quad \text{for }i=2,3,\dots 
\end{equation*}
Consider the $n-$th and the $(n+1)-$th outer iterations which contain $I_{1}$
and $I_{2}$ inner iterations, respectively. We define the sequence of
relative errors associated to these two outer iterations as 
\begin{equation}
e_{n,2},\dots ,e_{n,I_{1}},\widetilde{e}_{n+1,1},e_{n+1,2},\dots
,e_{n+1,I_{2}},  \label{eq:sequence}
\end{equation}%
where 
\begin{equation*}
\widetilde{e}_{n+1,1}=\frac{\Vert c_{n+1,1}-c_{n,I_{1}}\Vert _{L_{2}(\Omega
)}}{\Vert c_{n,I_{1}}\Vert _{L_{2}(\Omega )}}.
\end{equation*}

The inner iterations with respect to $i$ of the $n-$th outer iteration in
the Algorithm~\ref{alg:globalconv} are stopped if either $e_{n,2}<10^{-6}$
or $i=3$. Note that this rule is similar to the one used for
\textquotedblleft Test 2" in~\cite{Thanh2014}, where the maximal number of
inner iterations is set to be 5. We have observed in our numerical
experiments that the reconstruction results are essentially the same when we
use either 3 or 5 for the maximal number of inner iterations.

Concerning the outer iterations with respect to $n$ in Algorithm~\ref%
{alg:globalconv}, it can be seen from the stopping rule for the inner
iterations that each outer iteration consists of at least 2 and at most 3
inner iterations. Equivalently, the error sequence~\eqref{eq:sequence} has
at least 3 and at most 5 elements. We stop the outer iterations if there are
two consecutive outer iterations for which their error sequence~%
\eqref{eq:sequence} has three consecutive elements less than or equal to $%
5\times 10^{-4}$. 
We emphasize again that we have no rigorous justification for these stopping
rules; they rely on the content of the convergence theorem~\cite{Kliba2016}
and trial-and-error testing for simulated data.

We choose the final result for $c(\mathbf{x})$ by taking the average of its
approximations $c_{n,i}(\mathbf{x})$ corresponding to the relative errors in~%
\eqref{eq:sequence} that meet the stopping criterion for outer iterations.
The computed dielectric constant is determined as the maximal value of the
computed function $c(\mathbf{x})$. We have observed in our numerical studies
that we need no more than five outer iterations to obtain the final result.

\subsection{Reconstruction results}

In Table~\ref{tab:table2} we present the reconstruction results for the
dielectric constants of the scattering objects listed in Table~\ref%
{tab:table1}. Dielectric constants of targets were independently directly
measured by physicists from the Optoelectronics Center of the University of
North Carolina at Charlotte (M. A. Fiddy and S. Kitchin). The error in the
measurement are given by the standard deviation. The measured dielectric
constants were dependent on the frequency. Thus, in the third column of
Table~\ref{tab:table2} we present those directly measured dielectric
constants of the scattering objects at the corresponding optimal frequencies
(in column 2). Optimal frequencies are defined in section~\ref%
{sect:centralFreq}.

Our computed dielectric constants as well as the relative errors compared
with the measured dielectric constants (column 3) are displayed in the last
two columns of Table~\ref{tab:table2}. All relative errors vary
 between 1.02\% and 9.33\%. This is certainly a very good accuracy of
the reconstruction, especially given a very significant amount of noise in
the data (Section 5.2). We point out that a similar reconstruction accuracy has been
consistently demonstrated in all other works where performances of globally convergent
inversion methods were tested on experimental data~\cite{Klibanov:ip2010, Beili2012, Kuzhu2012, Thanh2014, Exp1, Exp2}. 
\begin{table}[h]
\caption{Measured and computed dielectric constant of the targets}
\label{tab:table2}\centering
\begin{tabular}{|c|c|c|c|c|}
\hline
Obj.\,ID & Freq. (GHz) & Measured $c$\,(std.~dev.) & Computed $c$ & Relative error   \\ \hline
1 & 3.10 & 4.50 (5.99\%) & 4.92 & 9.33\%   \\ \hline
2 & 3.01 & 5.45 (1.13\%) & 5.25 & 3.67\% \\ \hline
3 & 3.01 & 5.61 (21.3\%) & 5.12 & 8.73\%   \\ \hline
4 & 3.10 & 4.89 (2.89\%) & 4.94 & 1.02\%  \\ \hline
5 & 2.62 & 7.58 (4.69\%) & 8.18 & 7.92\% \\ \hline
6 & 2.62 & 4.80 (1.54\%) & 4.97 & 3.54\%   \\ \hline
\end{tabular}%
\end{table}

\begin{figure}[!h]
\centering
\subfloat[Projection of the (visualized) exact coefficient for object 6 on $\{y=0\}$ ]{\includegraphics[width=0.4\textwidth]{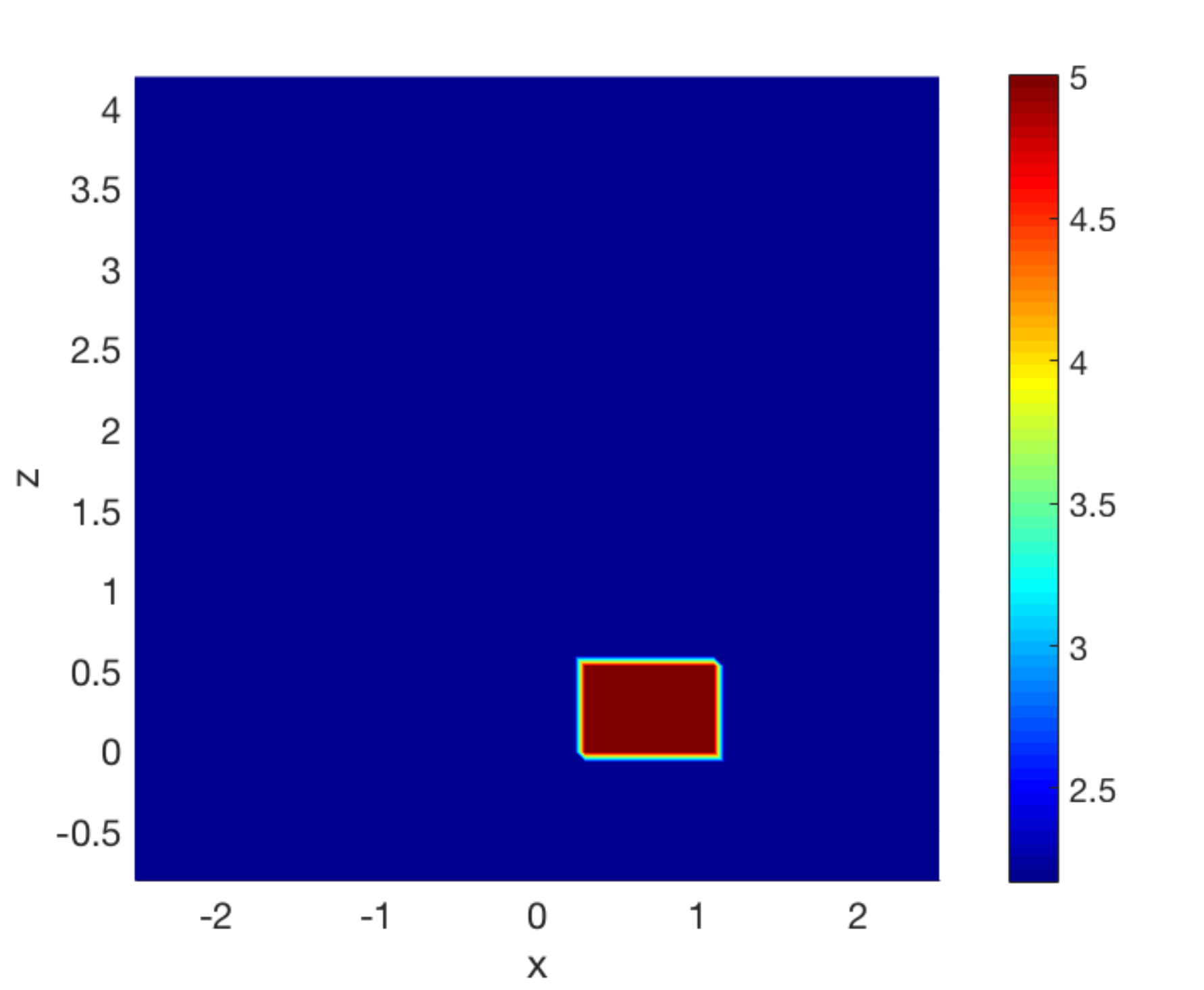}} \hspace{
0.2cm} 
\subfloat[Reconstruction  result of the coefficient for object 6, projected on $\{y=0\}$ ]{\includegraphics[width=0.4\textwidth]{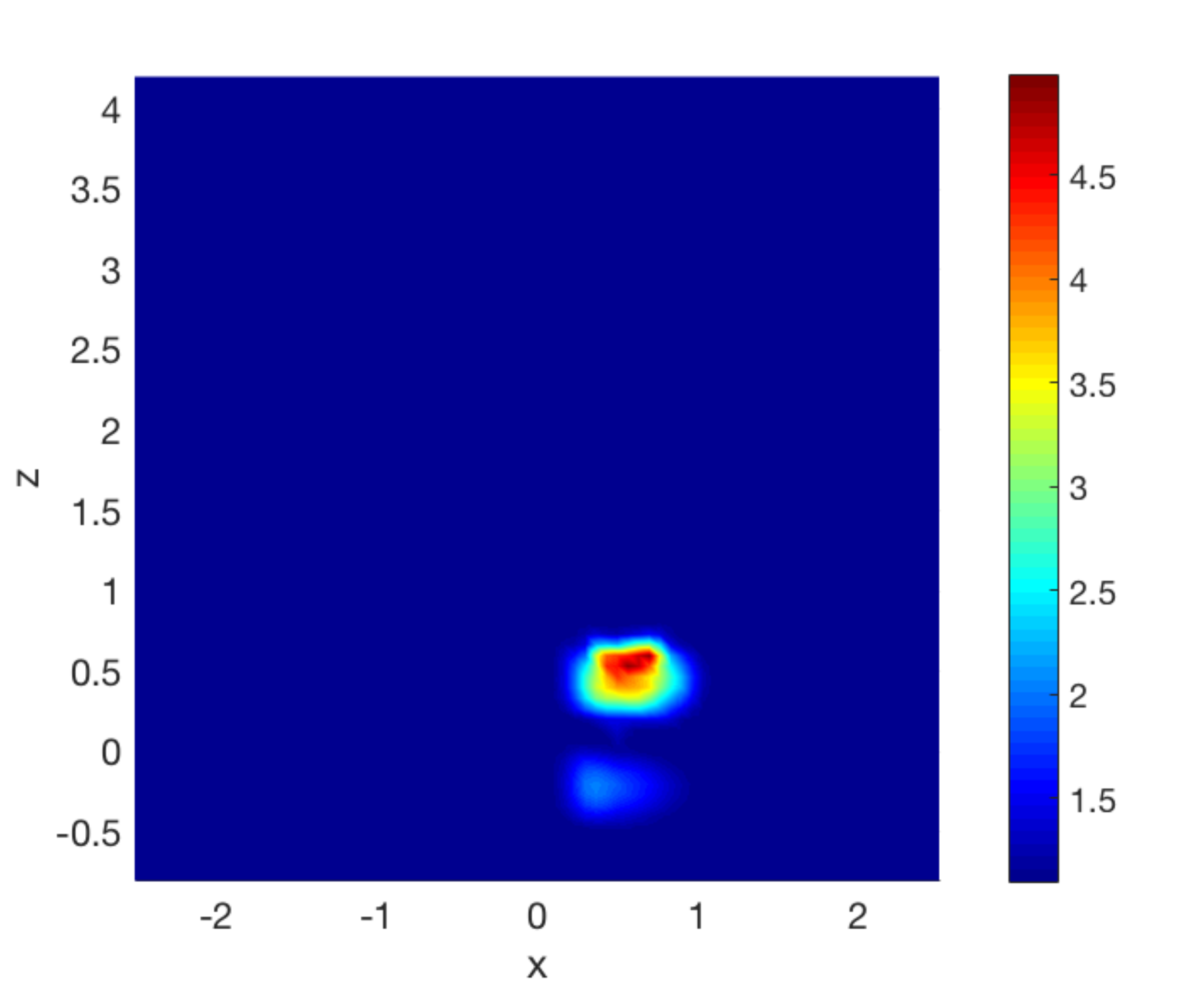}} \vspace{-0.0cm} 
\subfloat[Exact geometry of object 6]{\includegraphics[width=0.4\textwidth]{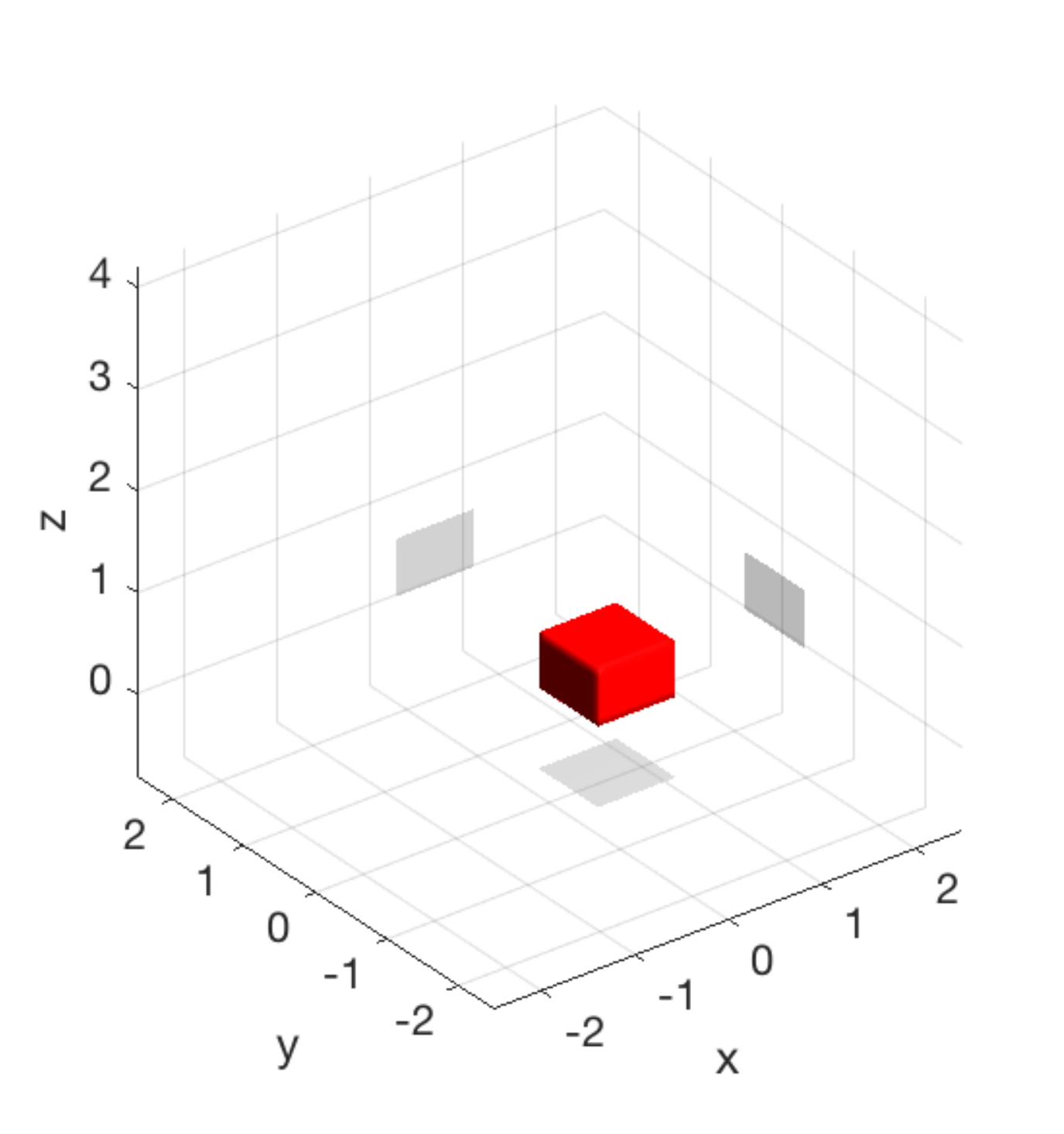}} \hspace{0.2cm} 
\subfloat[Reconstruction  result  of object 6]{\includegraphics[width=0.4\textwidth]{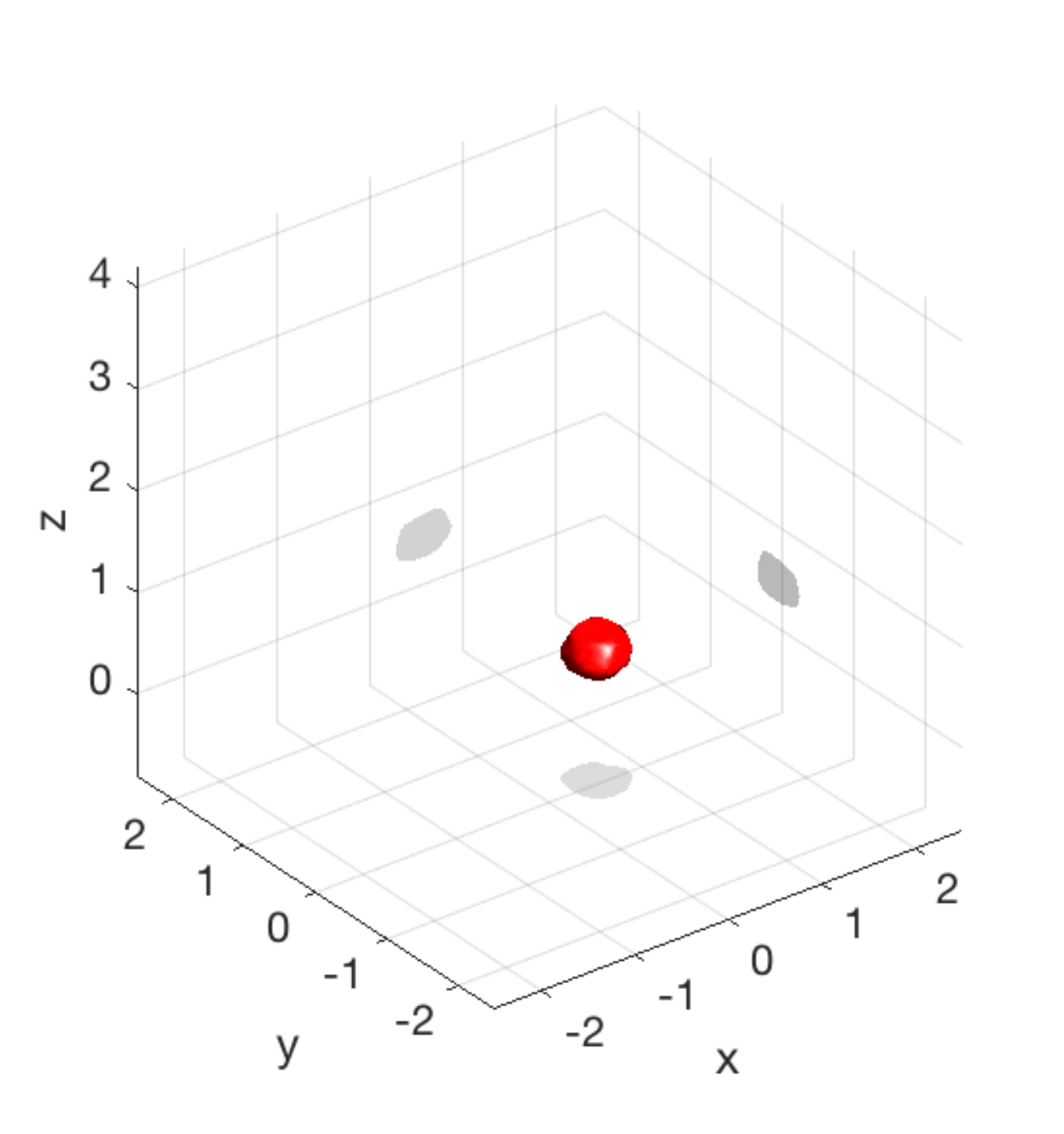}}
\caption{2D and 3D visualizations of exact (left) and reconstructed (right)
geometry for object 6.}
\label{fi:6}
\end{figure}

\begin{figure}[!h]
\centering
\subfloat[Exact geometry of object
2]{\includegraphics[width=0.4\textwidth]{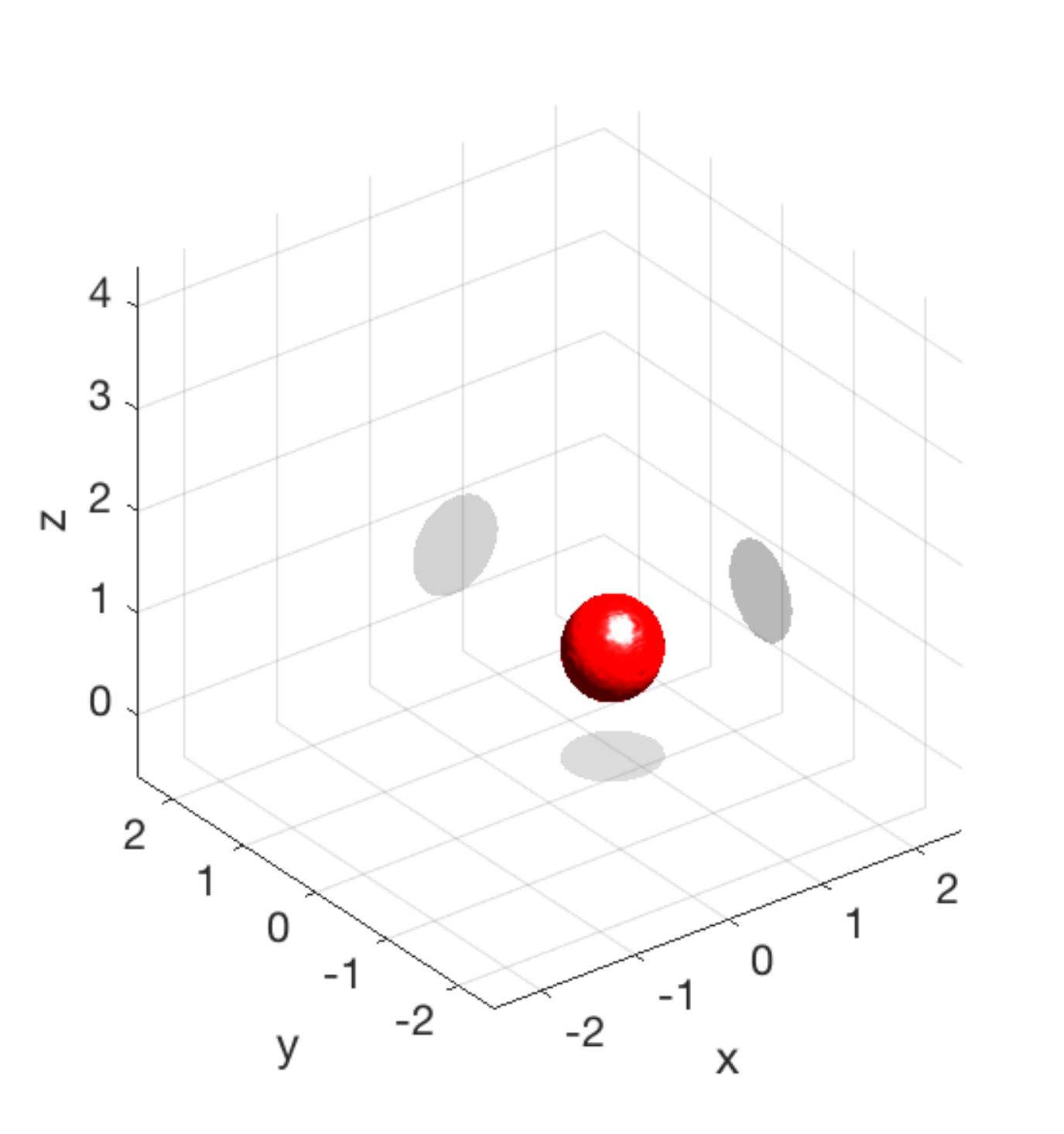}} \hspace{0.2cm} 
\subfloat[Reconstruction  result of object
2]{\includegraphics[width=0.4\textwidth]{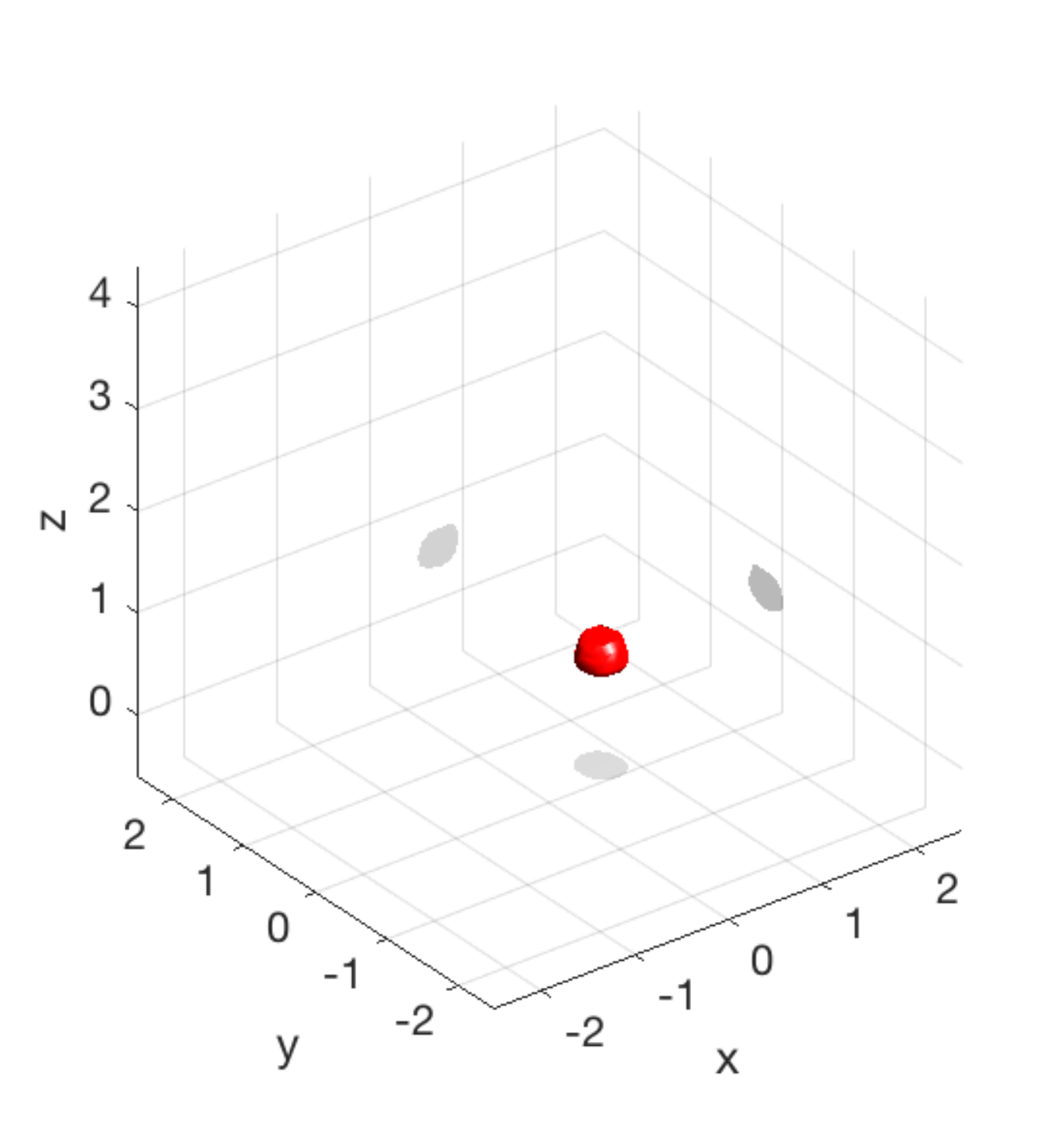}} \vspace{-0.0cm} 
\subfloat[Exact gometry of object
4]{\includegraphics[width=0.4\textwidth]{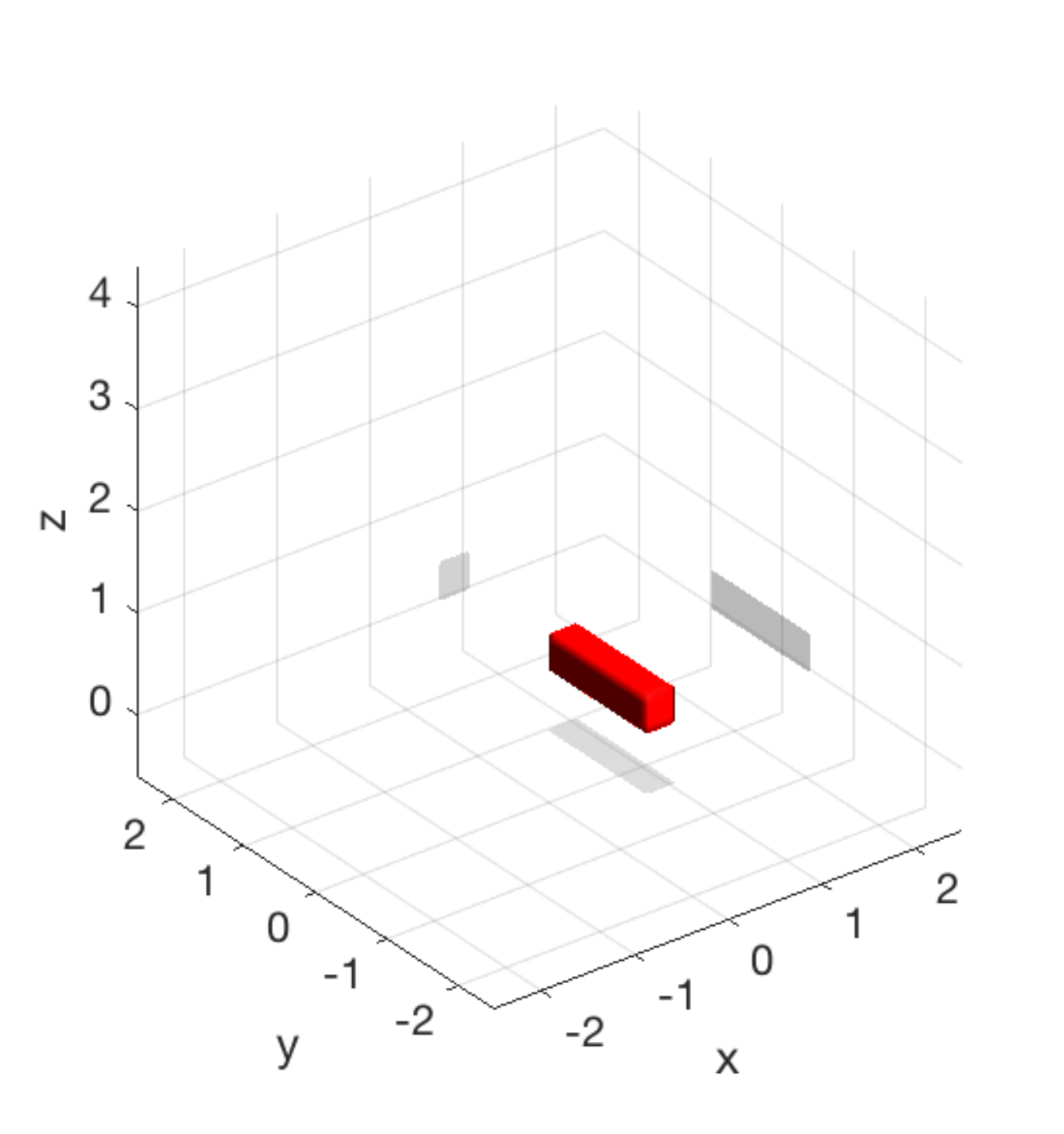}} \hspace{0.2cm} 
\subfloat[Reconstruction result of object
4]{\includegraphics[width=0.4\textwidth]{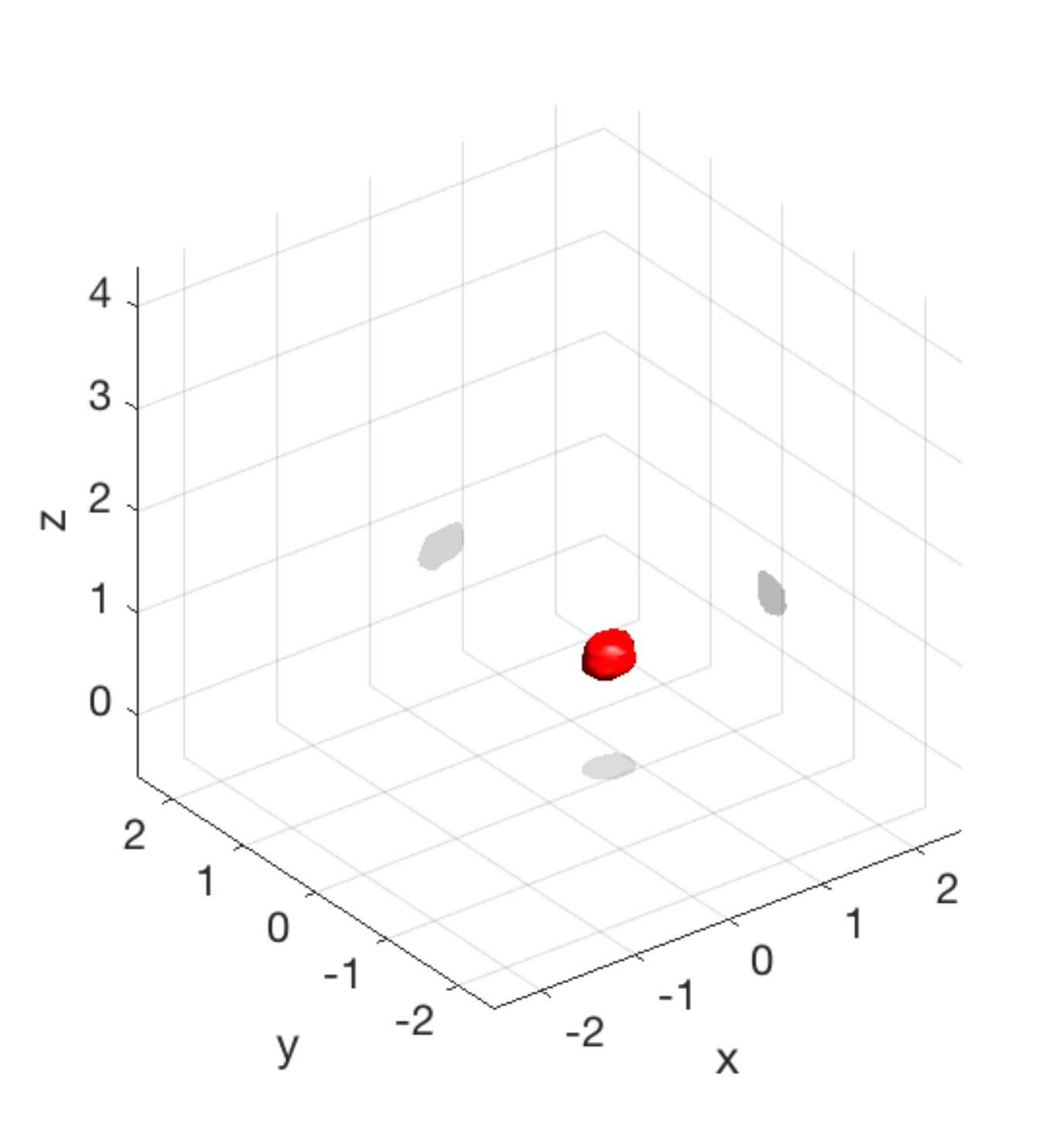}}
\caption{ 3D visualizations of exact (left) and reconstructed (right)
geometry for objects 2 and 4.}
\label{fi:7}
\end{figure}

Figures~\ref{fi:6} and \ref{fi:7} present the visualizations of exact and
reconstructed geometry of three out of six representative targets
using \texttt{isosurface} in MATLAB. We do not display the reconstructed
image for three other objects since they look similar to the ones
presented. For the isosurface plotting, we chose the isovalue as 50\% of the
maximal value of the computed function $c(\mathbf{x}) $. This
choice is based on our successful  experience in previous papers~\cite{Kliba2016, Exp1} 
and is applied to all other objects.
One can see from these figures that locations of the targets are well
reconstructed.

\section*{Acknowledgements}

The work of M. V. Klibanov and D.-L. Nguyen was supported by the US Army Research
Laboratory and US Army Research Office grant W911NF-15-1-0233 as well as by
the Office of Naval Research grant N00014-15-1-2330. The work of L. H. Nguyen
was partially supported by research funds provided by University of North
Carolina at Charlotte.

\bibliographystyle{plain}
\bibliography{ip-biblio}

\end{document}